# Small-time behavior of beta coalescents

Julien Berestycki[a], Nathanaël Berestycki[b] and Jason Schweinsberg[c,*]

[a]*Laboratoire d'Analyse, Topologie, Probabilités UMR 6632, Centre de Mathématiques et Informatique, Université de Provence, 39 rue F. Joliot-Curie, 13453 Marseille cedex 13, France*
[b]*University of British Columbia, Room 121 – 1984, Mathematics Road, Vancouver, BC VCT 1Z2, Canada*
[c]*Department of Mathematics, 0112, 9500 Gilman Drive, La Jolla, CA 92093-0112, USA. E-mail: jschwein@math.ucsd.edu*

**Abstract.** For a finite measure $\Lambda$ on $[0,1]$, the $\Lambda$-coalescent is a coalescent process such that, whenever there are $b$ clusters, each $k$-tuple of clusters merges into one at rate $\int_0^1 x^{k-2}(1-x)^{b-k}\Lambda(\mathrm{d}x)$. It has recently been shown that if $1 < \alpha < 2$, the $\Lambda$-coalescent in which $\Lambda$ is the Beta$(2-\alpha, \alpha)$ distribution can be used to describe the genealogy of a continuous-state branching process (CSBP) with an $\alpha$-stable branching mechanism. Here we use facts about CSBPs to establish new results about the small-time asymptotics of beta coalescents. We prove an a.s. limit theorem for the number of blocks at small times, and we establish results about the sizes of the blocks. We also calculate the Hausdorff and packing dimensions of a metric space associated with the beta coalescents, and we find the sum of the lengths of the branches in the coalescent tree, both of which are determined by the behavior of coalescents at small times. We extend most of these results to other $\Lambda$-coalescents for which $\Lambda$ has the same asymptotic behavior near zero as the Beta$(2-\alpha, \alpha)$ distribution. This work complements recent work of Bertoin and Le Gall, who also used CSBPs to study small-time properties of $\Lambda$-coalescents.

**Résumé.** L'objet de ce travail est l'étude du comportement asymptotique en temps petit des Beta-coalescents. Ces processus décrivent la limite d'échelle de la généalogie d'un certain nombre de modèles en génétique des populations. Nous donnons en particulier un théorème de convergence presque sûre pour le nombre de blocs renormalisé. Nous décrivons également le comportement asymptotique des tailles des blocs. Ces résultats permettent de calculer la dimension de Hausdorff et la dimension de packing d'un espace métrique associé à ce type de coalescents, ainsi que la longueur totale des branches de l'arbre de coalescence. Ce dernier résultat correspond à une question qui se pose en génétique des populations. Enfin, ces résultats sont en partie étendus par des arguments de couplage aux cas de $\Lambda$-coalescents pour lesquels la mesure $\Lambda$ a un comportement près de 0 semblable à celui d'une distribution Beta. Les méthodes employées reposent essentiellement sur un lien entre Beta-coalescent et les processus de branchement à espace d'état continu.

---

*Supported in part by NSF Grant DMS-05-04882.

## 1. Introduction

Coalescent processes are stochastic models of a system of particles that start out separated and then merge into clusters as time goes forward. Coalescent processes have applications to areas such as physical chemistry, where one can think of the merging of physical particles, astronomy, where we have the merging of galaxies into clusters, and biology, where ancestral lines of a sample from a population merge as we go backward in time. See [1, 5] for surveys.

Much work on coalescence has focused on processes in which only two clusters can merge at a time. However, Pitman [30] and Sagitov [31] introduced coalescents with multiple collisions, in which many clusters can merge at once into a single cluster. To define these processes precisely, let $\mathcal{P}_n$ be the set of partitions of $\{1, \ldots, n\}$, and let $\mathcal{P}$ be the set of partitions of $\mathbb{N}$. For all partitions $\pi \in \mathcal{P}$, let $R_n \pi$ be the restriction of $\pi$ to $\{1, \ldots, n\}$, meaning that $R_n \pi \in \mathcal{P}_n$, and two integers $i$ and $j$ are in the same block of $R_n \pi$ if and only if they are in the same block of $\pi$. A coalescent with multiple collisions is a $\mathcal{P}$-valued Markov process $(\Pi(t), t \geq 0)$ such that $\Pi(0)$ is the partition of $\mathbb{N}$ into singletons and, for all $n \in \mathbb{N}$, the process $(R_n \Pi(t), t \geq 0)$ is a $\mathcal{P}_n$-valued Markov process with the property that whenever there are $b$ blocks, each transition that involves merging $k$ blocks of the partition into one happens at rate $\lambda_{b,k}$, and these are the only possible transitions. The rates $\lambda_{b,k}$ do not depend on $n$ nor on the numbers of integers in the $b$ blocks. Pitman showed that the transition rates must satisfy

$$\lambda_{b,k} = \int_0^1 x^{k-2}(1-x)^{b-k} \Lambda(\mathrm{d}x) \qquad (1)$$

for some finite measure $\Lambda$ on $[0, 1]$, and the coalescent process such that (1) holds for a particular measure $\Lambda$ is called the $\Lambda$-coalescent. When $\Lambda$ is a unit mass at zero, then each transition involves the merger of exactly two blocks, and each such transition occurs at rate 1. This process is known as Kingman's coalescent and was introduced in [25].

There has been a considerable amount of work concerning applications of these processes. Sagitov [31] showed that coalescents with multiple collisions can describe the genealogy of populations in which there are occasionally very large families. See [29] for further results in this direction. Durrett and Schweinsberg [18] showed that coalescents with multiple collisions can be used to model the genealogy of a population that periodically experiences beneficial mutations. Schweinsberg [33] considered the genealogy of supercritical Galton–Watson processes in which the probabability of having $k$ or more offspring decays like $Ck^{-\alpha}$ for some constant $C$. When $1 \leq \alpha < 2$, the genealogy of this process, as the population size tends to infinity, converges to the $\Lambda$-coalescent in which $\Lambda$ is the Beta$(2 - \alpha, \alpha)$ distribution. Birkner et al. [11] established a continuous version of these results, showing that the $\Lambda$-coalescents that describe the genealogy of a continuous-state branching process (CSBP) are precisely those in which $\Lambda$ is the Beta$(2-\alpha, \alpha)$ distribution, where $0 < \alpha < 2$. The $\alpha = 1$ case had previously been established by Bertoin and Le Gall [6].

These results suggest that the $\Lambda$-coalescents in which $\Lambda$ is the Beta$(2 - \alpha, \alpha)$ distribution form an important one-parameter family of coalescents with multiple collisions that is worthy of further study. These results also suggest that it should be possible to use results about continuous-state branching processes to get new insight into the behavior of coalescent processes. The goal of this paper is to establish some results about the asymptotics of the Beta$(2 - \alpha, \alpha)$-coalescents at small times. Because the small-time behavior of $\Lambda$-coalescents depends only on properties of $\Lambda$ near zero, some of our results extend easily to $\Lambda$-coalescents that have the same behavior near zero as the Beta$(2 - \alpha, \alpha)$-coalescents, and we prove these results in this more general form. Note that when $\alpha = 1$, the Beta$(1, 1)$ distribution is the uniform distribution on $[0, 1]$, and the associated coalescent process, called the Bolthausen–Sznitman coalescent, has already been studied extensively (see, for example, [2, 9, 12, 22, 30]). We focus here on the case in which $1 < \alpha < 2$. Some of our results are closely related to results of Bertoin and Le Gall [8], who also used CSBPs to study the small-time behavior of $\Lambda$-coalescents.

### 1.1. Number of blocks

Our first result concerns the number of blocks at small times. We say the $\Lambda$-coalescent comes down from infinity if the number of blocks is a.s. finite for all $t > 0$, and stays infinite if the number of blocks is a.s.

infinite for all $t > 0$. It is known (see Example 15 of [32]) that the Beta$(2 - \alpha, \alpha)$ coalescent comes down from infinity if $1 < \alpha < 2$ and stays infinite if $0 < \alpha \leq 1$. Theorem 1.1 gives a limit theorem for the number of blocks at small times when $1 < \alpha < 2$. Note that this limit holds almost surely. Independently of the present work, Bertoin and LeGall obtained the limit in probability for a larger family of $\Lambda$-coalescents (see Lemma 3 of [8]).

**Theorem 1.1.** *Let $\Lambda$ be a finite measure on $[0,1]$ such that $\Lambda(\mathrm{d}x) = f(x)\,\mathrm{d}x$, where $f(x) \sim Ax^{1-\alpha}$ for some $\alpha \in (1,2)$ and $\sim$ means that the ratio of the two sides tends to one as $x \downarrow 0$. Let $(\Pi(t), t \geq 0)$ be the $\Lambda$-coalescent, and let $N(t)$ be the number of blocks of the partition $\Pi(t)$. Then*

$$\lim_{t \downarrow 0} t^{1/(\alpha-1)} N(t) = \left(\frac{\alpha}{A\Gamma(2-\alpha)}\right)^{1/(\alpha-1)} \quad a.s. \tag{2}$$

*In particular, if $\Lambda$ is the Beta$(2 - \alpha, \alpha)$ distribution, then*

$$\lim_{t \downarrow 0} t^{1/(\alpha-1)} N(t) = (\alpha \Gamma(\alpha))^{1/(\alpha-1)} \quad a.s. \tag{3}$$

To see how (3) follows from (2), note that for the Beta$(2 - \alpha, \alpha)$-coalescent, we have

$$\Lambda(\mathrm{d}x) = \frac{1}{\Gamma(2-\alpha)\Gamma(\alpha)} x^{1-\alpha}(1-x)^{\alpha-1}\,\mathrm{d}x,$$

so in this case $A = 1/[\Gamma(\alpha)\Gamma(2-\alpha)]$. Also, note that as $\alpha \uparrow 2$, the Beta$(2-\alpha, \alpha)$ distribution converges to the unit mass at zero. Consequently, although Theorem 1.1 is stated for $1 < \alpha < 2$, Kingman's coalescent can be viewed as corresponding to $\alpha = 2$. Indeed, it is known for Kingman's coalescent (see Section 4.2 of [1]) that $tN(t) \to 2$ a.s. as $t \downarrow 0$, which is what one gets plugging $\alpha = 2$ into (3).

In Section 4, Theorem 1.1 is obtained by relating the behavior to continuous-state branching processes. In [3], we present an alternative approach based on continuous stable random trees and the Kesten–Stigum theorem.

### 1.2. Block sizes

We now consider the sizes of the blocks of the beta coalescents. It is clear from the definition that if $(\Pi(t), t \geq 0)$ is a $\Lambda$-coalescent, then $\Pi(t)$ is an exchangeable random partition of $\mathbb{N}$ for all $t > 0$. It thus follows from results of Kingman [24] that if $B \subset \mathbb{N}$ is a block of the partition $\Pi(t)$, then the limit

$$\lim_{m \to \infty} \frac{1}{m} \sum_{i=1}^{m} \mathbf{1}_{\{i \in B\}}$$

exists almost surely and is called the asymptotic frequency of $B$. If, for each $t > 0$, the sum of the asymptotic frequencies of the blocks of $\Pi(t)$ equals one almost surely, then we say the coalescent has proper frequencies. Pitman showed (see Theorem 8 of [30]) that the $\Lambda$-coalescent has proper frequencies if and only if $\int_0^1 x^{-1}\Lambda(\mathrm{d}x) = \infty$. In particular, the Beta$(2-\alpha, \alpha)$-coalescent has proper frequencies if and only if $\alpha \geq 1$. When $0 < \alpha < 1$, for each $t > 0$, almost surely a positive asymptotic fraction of the integers will be in singleton blocks of $\Pi(t)$, so the sum of the asymptotic frequencies will be less than one. For coalescents with proper frequencies, almost surely $\Pi(t)$ has no singletons for all $t > 0$.

If $(\Pi(t), t \geq 0)$ is a $\Lambda$-coalescent, then one can construct a ranked $\Lambda$-coalescent $(\Theta(t), t \geq 0)$ such that $\Theta(t)$ is the sequence of asymptotic frequencies of the blocks of the partition $\Pi(t)$, ranked in decreasing order. For most $\Lambda$-coalescents, there appears to be no simple description of the distribution of $\Theta(t)$ for fixed $t$. An exception is the Bolthausen–Sznitman coalescent, in which case $\Theta(t)$ has the Poisson–Dirichlet distribution with parameters $(\mathrm{e}^{-t}, 0)$; see [12, 30], or see [22] for a short proof using recursive trees. Also, for Kingman's coalescent, if $T_k = \inf\{t: N(t) \leq k\}$ is the first time at which the coalescent has $k$ blocks, then the distribution of $\Theta(T_k)$ is uniform on the simplex $\Delta_k = \{(x_1, \ldots, x_k): x_1 \geq \cdots \geq x_k, x_1 + \cdots + x_k = 1\}$, as shown in [25]. Theorem 1.2 below gives a related result for the Beta$(2 - \alpha, \alpha)$-coalescents with $1 < \alpha < 2$. It is possible that these coalescent processes may never have exactly $k$ blocks, so we calculate the conditional distribution of $\Theta(T_k)$ given $N(T_k) = k$, which is the event that the coalescent has exactly $k$ blocks at some time.

**Theorem 1.2.** *Let $(\Pi(t), t \geq 0)$ be the Beta$(2 - \alpha, \alpha)$-coalescent, where $1 < \alpha < 2$, and let $(\Theta(t), t \geq 0)$ be the associated ranked coalescent. Let $N(t)$ be the number of blocks of $\Pi(t)$ at time $t$. Fix a positive integer $k$, and let $X'_1, \ldots, X'_k$ be i.i.d. random variables with distribution $\mu$, where the Laplace transform of $\mu$ is given by*

$$\int_0^\infty e^{-\lambda x} \mu(dx) = 1 - (1 + \lambda^{1-\alpha})^{-1/(\alpha-1)}. \tag{4}$$

*Let $X_1, \ldots, X_k$ be the values of $X'_1, \ldots, X'_k$ ranked in decreasing order. Let $S_k = X_1 + \cdots + X_k$. If $g: \Delta_k \to [0, \infty)$ is a nonnegative measurable function, then*

$$E[g(\Theta(T_k))|N(T_k) = k] = \frac{1}{E[S_k^{1-\alpha}]} E\left[S_k^{1-\alpha} g\left(\frac{X_1}{S_k}, \ldots, \frac{X_k}{S_k}\right)\right]. \tag{5}$$

To see how this result is related to the result for Kingman's coalescent, note that if $\alpha = 2$, the right-hand side of (4) becomes $1/(1 + \lambda)$, so $\mu$ is the exponential distribution with mean 1. If $X_1, \ldots, X_k$ are obtained by ranking $k$ i.i.d. random variables that have the exponential distribution with mean 1 and $S_k = X_1 + \cdots + X_k$, then $S_k$ is independent of $(X_1/S_k, \ldots, X_k/S_k)$. Consequently, the right-hand side of (5) becomes $E[g(X_1/S_k, \ldots, X_k/S_k)]$. Furthermore, the distribution of $(X_1/S_k, \ldots, X_k/S_k)$ is uniform on $\Delta_k$. Note also that $\Theta(T_k)$ is independent of $T_k$ for Kingman's coalescent because exactly two blocks coalesce during each merger, but this property does not hold for other $\Lambda$-coalescents.

*Remark 1.3.* *The distribution $\mu$ first arose in the work of Slack [35], where it was used to describe the family sizes of critical Galton–Watson processes with heavy-tailed offspring distribution, at large times when conditioned on survival. More precisely, recall that Yaglom's limit law [23, 36] states that for critical Galton–Watson processes with finite variance, the distribution of the number of offspring at time $n$, conditioned to be positive and then rescaled to have mean 1, converges to the exponential distribution with mean 1 as $n \to \infty$. When the offspring distribution is in the domain of attraction of a stable law of index $\alpha \in (1, 2)$ (and thus does not have finite variance), Slack showed that the distribution of the number of offspring in generation $n$, conditioned to be positive and then rescaled to have mean 1, converges to $\mu$ as $n \to \infty$, thus proving an analog of Yaglom's limit law for offspring distributions with infinite variance. The $\alpha$-stable CSBP for $\alpha \in (1, 2)$ arises as a limit of Galton–Watson processes whose offspring distribution is in the domain of attraction of a stable law (see [16, 26]). Since beta coalescents can be recovered from the genealogy of such continuous-state branching processes, it is natural that the same distribution $\mu$ arises here as well. Although our proof never uses it explicitly, many of our results can be understood intuitively in terms of Slack's theorem.*

We now consider the sizes of the blocks at small times. By evaluating the derivative of the right-hand side of (4) at zero, we see that $E[X'_i] = 1$ for all $i$. Therefore, $S_k$ will be approximately $k$ for large $k$. At small times, the number of blocks will be large, so Theorem 1.2 suggests that when there are $k$ blocks, the distribution of the asymptotic frequencies of these blocks will be approximately the distribution of $k$ independent random variables with distribution $\mu$, each divided by $k$. Theorem 1.4 below makes this observation rigorous. The motivation for this result comes from the recent work of Bertoin and Le Gall, who proved a similar statement (see Theorem 4 of [8]). Bertoin and Le Gall's result applies to a larger family of $\Lambda$-coalescents, as it requires only a regular variation condition on $\Lambda$ near zero. However, Bertoin and Le Gall prove only convergence in probability, whereas we establish almost sure convergence for the beta coalescents.

**Theorem 1.4.** *Let $(\Pi(t), t \geq 0)$ be the $\mathrm{Beta}(2-\alpha, \alpha)$-coalescent, where $1 < \alpha < 2$. Let $N(t,x)$ be the number of blocks of $\Pi(t)$ whose asymptotic frequency is at most $x$. Let $F(x) = \mu((0,x])$ for all $x$, where $\mu$ is the probability distribution defined in ([4](#)). Then*

$$\limsup_{t \downarrow 0} \sup_{x \geq 0} |t^{1/(\alpha-1)} N(t, t^{1/(\alpha-1)} x) - (\alpha \Gamma(\alpha))^{1/(\alpha-1)} F((\alpha \Gamma(\alpha))^{1/(\alpha-1)} x)| = 0 \quad a.s.$$

Note that by taking a limit as $x \to \infty$ in Theorem 1.4, we recover the result of Theorem 1.1 for the beta coalescents. Also, note that if $\alpha = 2$ and $\mu$ is the exponential distribution with mean 1, then the expression $(\alpha \Gamma(\alpha))^{1/(\alpha-1)} F((\alpha \Gamma(\alpha))^{1/(\alpha-1)} x)$ becomes $2(1 - \mathrm{e}^{-2x})$, and, as observed in [8], we again recover a known result for Kingman's coalescent (see Section 4.2 of [1]).

From Theorem 1.4, we obtain the following result for the size of the block containing the integer 1. Note that as a consequence of Kingman's work [24] on exchangeable random partitions, for coalescents with proper frequencies the asymptotic frequency of the block containing 1 is a size-biased pick from the asymptotic frequencies of all of the blocks.

**Proposition 1.5.** *Let $(\Pi(t), t \geq 0)$ be the $\mathrm{Beta}(2-\alpha, \alpha)$-coalescent, where $1 < \alpha < 2$, and let $K(t)$ be the asymptotic frequency of the block of $\Pi(t)$ containing 1. Then*

$$(\alpha \Gamma(\alpha))^{1/(\alpha-1)} t^{-1/(\alpha-1)} K(t) \xrightarrow{d} X \quad as\ t \downarrow 0,$$

*where $E[\mathrm{e}^{-\lambda X}] = (1 + \lambda^{\alpha-1})^{-\alpha/(\alpha-1)}$.*

Although the distribution $\mu$ has mean one and infinite variance, we can see by differentiating the Laplace transform of $X$ that $E[X] = \infty$. Also, as will be seen from the proof of the proposition, $X$ has the size-biased distribution $P(X \in \mathrm{d}x) = x\mu(\mathrm{d}x)$.

We also have the following result concerning the largest block of the coalescent at small times. While the size of a typical block and the block containing 1 are both of order $t^{1/(\alpha-1)}$, the size of the largest block is of order $t^{1/\alpha}$. This result follows from a Tauberian theorem, which gives information about the tail behavior of the distribution $\mu$, and extreme value theory. Recall that a random variable $X$ is said to have a Fréchet distribution of index $\alpha$ if $P(X \leq x) = \mathrm{e}^{-x^{-\alpha}}$ for all $x > 0$.

**Proposition 1.6.** *Let $(\Pi(t), t \geq 0)$ be the $\mathrm{Beta}(2-\alpha, \alpha)$-coalescent, where $1 < \alpha < 2$, and let $W(t)$ be the the largest of the asymptotic block frequencies of $\Pi(t)$. Then*

$$(\alpha \Gamma(\alpha) \Gamma(2-\alpha))^{1/\alpha} t^{-1/\alpha} W(t) \xrightarrow{d} X \quad as\ t \downarrow 0,$$

*where $X$ has the Fréchet distribution of index $\alpha$.*

This result suggests that there should be a whole range of block sizes between the typical size $t^{1/(\alpha-1)}$ and the largest block size $t^{1/\alpha}$. This is made more precise in [3] where we analyze the precise multifractal nature of the Beta-coalescent.

### 1.3. Hausdorff and packing dimensions

Given a $\mathcal{P}$-valued coalescent process, we can define a metric $d$ on $\mathbb{N}$ such that

$$d(i, j) = \inf\{t \colon i \text{ and } j \text{ are in the same block at time } t\}.$$

For all $i, j, k \in \mathbb{N}$, we have $d(i, j) \leq \max\{d(i, k), d(k, j)\}$, so $d$ is an ultrametric on $\mathbb{N}$. Let $(S, d)$ be the completion of $(\mathbb{N}, d)$, and note that the extension of $d$ to $S$ is also an ultrametric.

We now review the definitions of the Hausdorff and packing dimensions, following closely the discussions in [20, 21]. Let $(X, d)$ be a metric space. For $U \subset X$, let $|U| = \sup\{d(x, y) \colon x, y \in U\}$ denote the diameter of $U$. If $\{V_i\}_{i=1}^{\infty}$ is a collection of Borel sets such that $U \subset \bigcup_{i=1}^{\infty} V_i$, then we call $\{V_i\}_{i=1}^{\infty}$ a cover of $U$. If in addition $|V_i| \leq \delta$ for all $i$, then we call $\{V_i\}_{i=1}^{\infty}$ a $\delta$-cover of $U$. Given $s > 0$, the $s$-dimensional Hausdorff measure of $U$ is

$$\mathcal{H}_s(U) = \lim_{\delta \downarrow 0} \left( \inf \left\{ \sum_{i=1}^{\infty} |V_i|^s \colon \{V_i\}_{i=1}^{\infty} \text{ is a } \delta\text{-cover of } U \right\} \right).$$

The Hausdorff dimension of $U$ is

$$\dim_H(U) = \inf\{s \geq 0 \colon \mathcal{H}_s(U) = 0\} = \sup\{s \geq 0 \colon \mathcal{H}_s(U) = \infty\}.$$

If $\{V_i\}_{i=1}^{\infty}$ is a collection of disjoint open balls centered in $U$ such that $|V_i| \leq \delta$ for all $i$, then we say $\{V_i\}_{i=1}^{\infty}$ a $\delta$-packing of $U$. The $s$-dimensional packing premeasure of $U$ is

$$P_s(U) = \lim_{\delta \downarrow 0} \left( \sup \left\{ \sum_{i=1}^{\infty} |V_i|^s \colon \{V_i\}_{i=1}^{\infty} \text{ is a } \delta\text{-packing of } U \right\} \right). \tag{6}$$

The $s$-dimensional packing outer measure of $U$ is then defined to be

$$p_s(U) = \inf \left\{ \sum_{i=1}^{\infty} P_s(V_i) \colon \{V_i\}_{i=1}^{\infty} \text{ is a cover of } U \right\}. \tag{7}$$

The packing dimension of $U$ is

$$\dim_P(U) = \inf\{s \geq 0 \colon p_s(U) = 0\} = \sup\{s \geq 0 \colon p_s(U) = \infty\}. \tag{8}$$

The Hausdorff dimension of a set is always less than or equal to the packing dimension (see, for example, Chapter 3 of [21]).

Evans [20] investigated the fractal properties of the metric space associated with Kingman's coalescent. He showed that the Hausdorff and packing dimensions are both equal to one almost surely, and that the metric space is capacity equivalent to the unit interval. Donnelly et al. [13] showed that, for a coalescent process resulting from coalescing Brownian motions on the circle, the associated completed metric space $(S, d)$ almost surely has Hausdorff and packing dimensions of $1/2$ and is capacity-equivalent to the middle-$\frac{1}{2}$ Cantor set. Our next result implies that the Hausdorff and packing dimensions of the metric space associated with the Beta$(2 - \alpha, \alpha)$ coalescent with $1 < \alpha < 2$ has Hausdorff and packing dimensions equal to $1/(\alpha - 1)$. Note that again we get the correct result for Kingman's coalescent by substituting $\alpha = 2$.

**Theorem 1.7.** *Let $\Lambda$ be a finite measure on $[0, 1]$ satisfying the conditions of Theorem 1.1. Let $(S, d)$ be the metric space associated with the $\Lambda$-coalescent $(\Pi(t), t \geq 0)$. Then, the Hausdorff and packing dimensions of $S$ are both $1/(\alpha - 1)$ almost surely.*

### 1.4. Dynamics of the number of blocks

Theorem 1.1 gives an almost sure limit theorem for the number of blocks in the coalescent at small times. Here we consider in more detail the dynamics of the process $(N(t), t \geq 0)$, for $\Lambda$-coalescents satisfying the assumptions of Theorem 1.1.

Let $\zeta_{n,k}$ be the probability that, if the $\Lambda$-coalescent has $n$ blocks, then it will lose exactly $k$ blocks at the time of the next merger. More precisely, let $(\Pi_n(t), t \geq 0)$ be the $\Lambda$-coalescent restricted to $\{1, \ldots, n\}$, and let $N_n(t)$ be the number of blocks of $\Pi_n(t)$. If $T = \inf\{t \colon \Pi_n(t) \neq \Pi_n(0)\}$, then $\zeta_{n,k} = P(N_n(T) = n - k)$. Note that if $\lambda_{n,k}$ is given by (1) and $\lambda_n = \sum_{k=2}^{n} \binom{n}{k} \lambda_{n,k}$ is the total merger rate when the coalescent has $n$ blocks, then

$$\zeta_{n,k} = \binom{n}{k+1} \frac{\lambda_{n,k+1}}{\lambda_n},$$

because $k+1$ blocks have to merge for the number of blocks to be reduced by $k$. It is not difficult to calculate (see Lemma 4 in [8]) that

$$\lim_{n\to\infty} \zeta_{n,k} = \frac{\alpha\Gamma(k+1-\alpha)}{(k+1)!\Gamma(2-\alpha)}. \tag{9}$$

If we define $\zeta_k = \alpha\Gamma(k+1-\alpha)/[(k+1)!\Gamma(2-\alpha)]$, then $\sum_{k=1}^{\infty} \zeta_k = 1$ and $\sum_{k=1}^{\infty} k\zeta_k = 1/(\alpha-1)$, as shown in Eqs (39) and (40) of [8]. Therefore, there is a probability distribution, which we call $\zeta$, on the positive integers corresponding to $(\zeta_k)_{k=1}^{\infty}$, and this distribution has mean $1/(\alpha-1)$. Thus, at small times, when the number of blocks is large, the successive jumps of the process $(N(t), t \geq 0)$ are approximately i.i.d. with distribution $\zeta$. We can use a renewal argument to establish the following theorem.

**Theorem 1.8.** *Let $\Lambda$ be a finite measure on $[0,1]$ satisfying the conditions of Theorem 1.1. Let $(\Pi(t), t \geq 0)$ be the $\Lambda$-coalescent. Let $N(t)$ be the number of blocks of $\Pi(t)$, and let $V_n$ be the event that $N(t) = n$ for some $t$. Then*

$$\lim_{n\to\infty} P(V_n) = \alpha - 1.$$

Once again, the case $\alpha = 2$ corresponds to Kingman's coalescent, where $P(V_n) = 1$ for all $n$ because the process $(N(t), t \geq 0)$ visits every integer. As $\alpha$ gets smaller, there are more large mergers that cause $(N(t), t \geq 0)$ to skip over some integers.

*1.5. Total time in the tree*

Given a $\Lambda$-coalescent $(\Pi(t), t \geq 0)$, consider the process $(R_n\Pi(t), t \geq 0)$, which is the coalescent restricted to $\{1, \ldots, n\}$ so the process starts with just $n$ blocks. For $k = 2, \ldots, n$, let $D_k$ be the duration of time for which $\Pi(t)$ has exactly $k$ blocks. Then

$$L_n = \sum_{k=2}^{n} kD_k$$

is the sum of the lengths of all the branches in the coalescent tree. This quantity has biological significance because if the coalescent process represents the ancestral tree of a sample of $n$ individuals from the population and a mutation occurs along one of the branches of this tree, then the $n$ individuals in the sample will not all have the same gene at the site of the mutation. Consequently, if mutations occur at rate $\theta$ along each branch and each mutation happens at a different site, then the number of "segregating sites" at which the $n$ sampled individuals do not have the same gene should be approximately $\theta L_n$.

For Kingman's coalescent, it is easily verified that

$$\frac{L_n}{\log n} \xrightarrow{p} 2,$$

where $\to_p$ denotes convergence in probability. Durrett and Schweinsberg [18] studied the case in which $\Lambda$ has a unit mass at zero as well as a component that allows for multiple mergers. Möhle [28] obtained a recursive equation for the limiting distribution of $n^{-1}L_n$ under the condition $\int_0^1 x^{-2}\Lambda(\mathrm{d}x) < \infty$, which implies that the total merger rate is finite even when the number of blocks is infinite. The result below includes the Beta$(2-\alpha, \alpha)$ coalescents for $1 < \alpha < 2$.

**Theorem 1.9.** *Let $\Lambda$ be a finite measure on $[0,1]$ satisfying the conditions of Theorem 1.1. Let $L_n$ be as defined above for the $\Lambda$-coalescent. Then*

$$\frac{L_n}{n^{2-\alpha}} \xrightarrow{p} \frac{\alpha(\alpha-1)}{A\Gamma(2-\alpha)(2-\alpha)}.$$

In [3], more precise results on the structure of the population under such a model are obtained using an approach based on continuous random trees.

The rest of this paper is organized as follows. In Section 2, we review some facts about continuous-state branching processes that we will need, and state the connection between CSBPs and Beta$(2-\alpha,\alpha)$ coalescents that was established in [11]. In Section 3, we record some results that will allow us to couple two coalescents with multiple collisions, which will be used to extend some of our results beyond the beta coalescents. We prove Theorem 1.1 in Section 4. Theorems 1.2 and 1.4 and Propositions 1.5 and 1.6 are proved in Section 5. We prove Theorem 1.7 in Section 6 and Theorems 1.8 and 1.9 in Section 7.

## 2. Beta coalescents and continuous-state branching processes

In this section, we review the results in [11] that relate continuous-state branching processes to beta coalescents. Continuous-state branching processes are the continuous versions of Galton–Watson processes. More formally, a continuous-state branching process is a $[0,\infty]$-valued Markov process $(Z(t), t \geq 0)$ whose transition functions $p_t(x, \cdot)$ satisfy

$$p_t(x+y, \cdot) = p_t(x, \cdot) * p_t(y, \cdot) \quad \text{for all } x, y \geq 0. \tag{10}$$

That is, the sum of independent copies of the process started at $x$ and $y$ has the same distribution as the process started at $x + y$. We think of $Z(t)$ as being the size of a population at time $t$, and the property (10) is called the branching property, because it can loosely be interpreted as meaning that if we start with a population size of $x + y$, then number of offspring of the first $x$ individuals is independent of the number of offspring of the remaining $y$.

For each $t \geq 0$, there is a function $u_t : [0, \infty) \to \mathbb{R}$ such that

$$E[e^{-\lambda Z(t)} | Z_0 = a] = e^{-au_t(\lambda)}. \tag{11}$$

If we exclude processes with an instantaneous jump to infinity, the functions $u_t$ satisfy the differential equation

$$\frac{\partial u_t(\lambda)}{\partial t} = -\Psi(u_t(\lambda)), \tag{12}$$

where $\Psi : [0, \infty) \to \mathbb{R}$ is a function of the form

$$\Psi(u) = \alpha u + \beta u^2 + \int_0^\infty (e^{-xu} - 1 + xu\mathbf{1}_{\{x \leq 1\}})\pi(\mathrm{d}x), \tag{13}$$

where $\alpha \in \mathbb{R}$, $\beta \geq 0$, and $\pi$ is a Lévy measure on $(0, \infty)$ satisfying $\int_0^\infty (1 \wedge x^2)\pi(\mathrm{d}x) < \infty$. The function $\Psi$ is called the branching mechanism of the CSBP.

As shown in [6], one can extend the CSBP to a two-parameter process $(Z(t,a), t \geq 0, a \geq 0)$ such that $Z(0, a) = a$ for all $a \geq 0$ and, for all $a, b \geq 0$, the process $(Z(t, a+b) - Z(t,a), t \geq 0)$ is independent of $(Z(t,c), t \geq 0, 0 \leq c \leq a)$ and has the same law as a CSBP with branching mechanism $\Psi$ started at $b$. Here, we think of $Z(t,a)$ as the number of individuals at time $t$ descended from the first $a$ individuals at time zero. For fixed $t$, the process $(Z(t,a), a \geq 0)$ is a subordinator, and it then follows from (11) that the Laplace exponent of this subordinator is the function $\lambda \mapsto u_t(\lambda)$.

Along the same lines, one can work with a measure-valued process $(M_t, t \geq 0)$ taking its values in the set of finite measures on $[0, 1]$ such that $(M_t([0, a]), t \geq 0, 0 \leq a \leq 1)$ has the same finite-dimensional distributions as $(Z(t,a), t \geq 0, 0 \leq a \leq 1)$. Now $(M_t([0, a]), 0 \leq a \leq 1)$ is a subordinator with Laplace exponent $\lambda \mapsto u_t(\lambda)$ run for time 1, and if we set $Z(t) = M_t([0, 1])$, then $(Z(t), t \geq 0)$ is a CSBP with branching mechanism $\Psi$ started at 1. An explicit construction of $(M_t, t \geq 0)$ can be given using the lookdown construction of Donnelly and Kurtz [14]. See also Section 2 of [11] for a review of this construction in the $\beta = 0$ case. In [3], a construction of this process is obtained in terms of continuous stable random trees.

For the purposes of studying Beta$(2-\alpha,\alpha)$-coalescents with $1 < \alpha < 2$, we will consider CSBPs that have the stable branching mechanism $\Psi(\lambda) = \lambda^\alpha$, where $1 < \alpha < 2$. In this case, the Lévy measure is given by

$$\pi(\mathrm{d}x) = \frac{\alpha(\alpha-1)}{\Gamma(2-\alpha)} x^{-1-\alpha}\,\mathrm{d}x$$

(see, for example, Example 4 of [19]). Birkner et al. [11] showed that after a time change, the genealogy of this CSBP can be described by the Beta$(2-\alpha,\alpha)$-coalescent. The full construction of the beta coalescent relies on the lookdown construction of Donnelly and Kurtz [14]. We describe here an identity involving one-dimensional distributions which will be sufficient for the applications in this paper. We assume that $Z(t) = M_t([0,1])$, where $(M_t, t \geq 0)$ is the measure-valued process defined above. To define the time change, for all $t \geq 0$ let

$$R(t) = \alpha(\alpha-1)\Gamma(\alpha) \int_0^t Z(s)^{1-\alpha}\,\mathrm{d}s,$$

and let $R^{-1}(t) = \inf\{s\colon R(s) > t\}$. Note that in [11], the time change is given only up to a constant in Theorem 1.1, but one can determine the exact constant, for example, from the proof of Lemma 3.7 in [11]. Theorem 1.1 of [11] states that the process $(M_{R^{-1}(t)}/Z(R^{-1}(t)), t \geq 0)$ has the same law as the $\Lambda$-Fleming–Viot process introduced in [7], where $\Lambda$ is the Beta$(2-\alpha,\alpha)$ distribution. The following lemma then follows immediately from the duality discovered in [7] between the $\Lambda$-Fleming–Viot process and the $\Lambda$-coalescent.

**Lemma 2.1.** *If $(\Pi(t), t \geq 0)$ is a* Beta$(2-\alpha,\alpha)$*-coalescent and $(\Theta(t), t \geq 0)$ is the associated ranked coalescent, then for all $t > 0$, the distribution of $\Theta(t)$ is the same as the distribution of the sizes of the atoms of the measure $M_{R^{-1}(t)}/Z(R^{-1}(t))$, ranked in decreasing order.*

Lemma 2.2 below describes the number and sizes of the atoms of $M_t$, when the CSBP has the stable branching mechanism $\Psi(\lambda) = \lambda^\alpha$. This result, in combination with Lemma 2.1, will be the key to using continuous-state branching processes to get information about the number and sizes of the blocks of beta coalescents. Note that, as will be seen from the proof, the Lévy measure of the subordinator $\lambda \mapsto u_t(\lambda)$ is finite for all $t > 0$, so $M_t$ has only finitely many atoms.

**Lemma 2.2.** *Assume $\Psi(\lambda) = \lambda^\alpha$. Let $D(t)$ be the number of atoms of $M_t$, and let $J(t) = (J_1(t), \ldots, J_{D(t)}(t))$ be the sizes of the atoms of $M_t$, ranked in decreasing order. Then $D(t)$ is Poisson with mean $\theta_t = [(\alpha-1)t]^{-1/(\alpha-1)}$. Conditional on $D(t) = k$, the distribution of $J(t)$ is the same as the distribution of $(\theta_t^{-1}X_1, \ldots, \theta_t^{-1}X_k)$, where $X_1, \ldots, X_k$ are obtained by picking $k$ i.i.d. random variables with distribution $\mu$, and then ranking them in decreasing order.*

**Proof.** When $\Psi(\lambda) = \lambda^\alpha$, it is possible to solve (12) explicitly with the initial condition $u_0(\lambda) = 1$, and we get

$$u_t(\lambda) = [(\alpha-1)t + \lambda^{1-\alpha}]^{-1/(\alpha-1)}$$

(see Eq. (2.15) of [26]). The sizes of the atoms of $M_t$ are precisely the sizes of the jumps of a subordinator with Laplace exponent $\lambda \mapsto u_t(\lambda)$ run for time 1. The number of atoms of $M_t$ is the number of jumps of this subordinator, which has the Poisson distribution with some mean $\theta_t$. Note that since $\lim_{\lambda \to \infty} \lambda^{-1} u_t(\lambda) = 0$, the subordinator has no drift (see the formula at the bottom of p. 72 in [4]). Therefore, $P(Z(t) = 0) = \mathrm{e}^{-\theta_t}$. Using (11) with $a = 1$ and the Monotone Convergence Theorem,

$$P(Z(t) = 0) = \lim_{\lambda \to \infty} E[\mathrm{e}^{-\lambda Z(t)}] = \lim_{\lambda \to \infty} \mathrm{e}^{-u_t(\lambda)}.$$

It follows that

$$\theta_t = \lim_{\lambda \to \infty} u_t(\lambda) = [(\alpha-1)t]^{-1/(\alpha-1)}. \tag{14}$$

Since $\theta_t < \infty$, the subordinator has jumps at only a finite rate, and the sizes of the jumps are nonnegative i.i.d. random variables, whose distributions can be read from the Laplace exponent of the subordinator.

To obtain the distribution of the jump sizes, let $\mu_t$ be the distribution of $\theta_t^{-1} X$, where $X$ has distribution $\mu$. Then

$$\int_0^\infty (1 - e^{-\lambda x}) \theta_t \mu_t(dx) = \theta_t \left(1 - \int_0^\infty e^{-\lambda x} \mu_t(dx)\right) = \theta_t \left(1 - \int_0^\infty e^{-\lambda \theta_t^{-1} x} \mu(dx)\right)$$
$$= \theta_t (1 + (\lambda \theta_t^{-1})^{1-\alpha})^{-1/(\alpha-1)} = u_t(\lambda).$$

It follows that $\theta_t \mu_t(dx)$ is the Lévy measure of the subordinator, and therefore the subordinator has jumps with size distribution $\mu_t$ at rate $\theta_t$. This implies the lemma. □

Note that since the number of atoms in $M_t$ represents the number of individuals at time zero that have descendants alive at time $t$, the number $D(t)$ of atoms of $M_t$ is almost surely a decreasing function of $t$. This is clear, for example, from the construction in [11]. Furthermore, as a consequence of the branching property of CSBPs, if

$$M_t = \sum_{i=1}^{D(t)} J_i(t) \delta_{a_i},$$

where $\delta_{a_i}$ denotes a unit mass at $a_i$, then conditional on $(M_s, 0 \leq s \leq t)$, the processes $(M_{t+s}(\{a_i\}), s \geq 0)$ for $i = 1, \ldots, D(t)$ have the same joint law as $D(t)$ independent CSBPs with branching mechanism $\Psi$ started from $J_1(t), \ldots, J_{D(t)}(t)$. Also, almost surely $M_{t+s}(\{a\}) = 0$ for all $s > 0$ and $a \notin \{a_1, \ldots, a_{D(t)}\}$.

Finally, we recall that every CSBP can be obtained as a time-change of a Lévy process with no negative jumps, as shown in [27, 34]. Given $\Psi$ as in (13), let $(Y(t), t \geq 0)$ be a Lévy process such that $Y(0) = a$ and $E[e^{-\lambda Y(t)}] = e^{-\lambda a + t\Psi(\lambda)}$. Define $(\tilde{Y}(t), t \geq 0)$ to be the process $(Y(t), t \geq 0)$ stopped when it hits zero. Let $U(t) = \inf\{s: \int_0^s \tilde{Y}(u)^{-1} du > t\}$. Then, if $(Z(t), t \geq 0)$ is a CSBP with branching mechanism $\Psi$ and $Z(0) = a$, the processes $(Z(t), t \geq 0)$ and $(Y(U(t)), t \geq 0)$ have the same law, if we adopt the convention that $Y(\infty) = \infty$.

## 3. Coupling of coalescent processes

To extend our results for the Beta$(2 - \alpha, \alpha)$-coalescents to other $\Lambda$-coalescents, it will be important to have techniques for coupling two coalescents with multiple collisions. To carry out this coupling, we will use the Poisson process construction of $\Lambda$-coalescents introduced by Pitman [30]. For simplicity, we assume that $\Lambda(\{0\}) = 0$, which will be the case in our examples.

Let $Q_x$ denote the distribution of an infinite sequence of $\{0, 1\}$-valued random variables that are one with probability $x$ and zero with probability $1 - x$. Let $L$ be the measure on $\{0, 1\}^\infty$ such that $L(B) = \int_0^1 Q_x(B) x^{-2} \Lambda(dx)$ for all measurable sets $B$. We will construct the $\Lambda$-coalescent from a Poisson point process on $[0, \infty) \times \{0, 1\}^\infty$ with intensity measure $dt \times L(d\xi)$. To do this, we first fix a positive integer $n$ and construct a $\mathcal{P}_n$-valued process $(\Pi_n(t), t \geq 0)$. We set $\Pi_n(0)$ to be the partition of $\{1, \ldots, n\}$ into singletons. If $(t, \xi)$ is a point of the Poisson process and $B_1, \ldots, B_b$ are the blocks of $\Pi_n(t-)$, ranked in order by their smallest element, then we define $\Pi_n(t)$ to be the partition obtained from $\Pi_n(t-)$ by merging all of the blocks $B_i$ such that $\xi_i = 1$, where we write $\xi = (\xi_1, \xi_2, \ldots)$. Since $\Lambda$ is a finite measure, it is easy to verify that for any fixed $t$, there are only finitely many points $(s, \xi)$ such that $s \leq t$ and at least two of $\xi_1, \ldots, \xi_n$ equal one. Consequently, the process $(\Pi_n(t), t \geq 0)$ is well defined. Furthermore, these processes are defined consistently for different values of $n$, which means there exists a unique $\mathcal{P}$-valued process $(\Pi(t), t \geq 0)$ such that $\Pi_n(t) = R_n \Pi(t)$ for all $n$ and $t$. The process $(\Pi(t), t \geq 0)$ is the $\Lambda$-coalescent, as shown in [30].

Below are the two coupling lemmas that we will use. Lemma 3.1 allows us to restrict our attention to the behavior of $\Lambda$ in a neighborhood of zero when we are concerned with small-time asymptotics of $\Lambda$-coalescents. This result appears implicitly in [32], but we give the short proof for completeness. Lemma 3.2 will allow us to compare other $\Lambda$-coalescents to beta coalescents.

**Lemma 3.1.** *Suppose $\Lambda_1$ and $\Lambda_2$ are finite measures on $[0,1]$ such that $\Lambda_1(\{0\}) = \Lambda_2(\{0\}) = 0$ and, for some $\delta > 0$, the restriction of $\Lambda_1$ to $[0, \delta]$ equals the restriction of $\Lambda_2$ to $[0, \delta]$. Then there exist $\mathcal{P}$-valued processes $(\Pi_1(t), t \geq 0)$ and $(\Pi_2(t), t \geq 0)$ such that $\Pi_1$ is a $\Lambda_1$-coalescent, $\Pi_2$ is a $\Lambda_2$-coalescent, and for some random time $t > 0$, we have $\Pi_1(s) = \Pi_2(s)$ for all $s < t$.*

**Proof.** For $i = 1, 2$, let $\Lambda_i'$ be the restriction of $\Lambda_i$ to $(\delta, 1]$, and let $\Lambda_3'$ be the restriction of $\Lambda_1$ to $[0, \delta]$. Let $\Psi_1'$, $\Psi_2'$, and $\Psi_3'$ be independent Poisson point processes on $[0, \infty) \times \{0,1\}^\infty$ such that $\Psi_i'$ has intensity $dt \times L_i(d\xi)$, where $L_i(B) = \int_0^1 Q_x(B) x^{-2} \Lambda_i'(dx)$ for all measurable $B$. Let $\Psi_1$ be the Poisson point process consisting of all points in $\Psi_1'$ and $\Psi_3'$, and let $\Psi_2$ be the Poisson point process consisting of all points in $\Psi_2'$ and $\Psi_3'$. For $i = 1, 2$, let $(\Pi_i(t), t \geq 0)$ be the $\mathcal{P}$-valued coalescent process obtained from $\Psi_i$ as described above. Then $(\Pi_1(t), t \geq 0)$ is a $\Lambda_1$-coalescent and $(\Pi_2(t), t \geq 0)$ is a $\Lambda_2$-coalescent.

For $i = 1, 2$, the total mass of $L_i$ is $\int_\delta^1 x^{-2} \Lambda_i(dx) \leq \delta^{-2} \Lambda_i([\delta, 1]) < \infty$. Thus, for $i = 1, 2$, if we define $t_i = \min\{s: (s, \xi) \text{ is a point of } \Psi_i'\}$, then $t_i > 0$. Therefore, if $t = \min\{t_1, t_2\}$, then the restrictions of $\Psi_1$ and $\Psi_2$ to $[0, t) \times \{0, 1\}^\infty$ are the same. It now follows from the construction that $\Pi_1(s) = \Pi_2(s)$ for all $s < t$. □

**Lemma 3.2.** *Suppose $\Lambda_1$ and $\Lambda_2$ are finite measures on $[0,1]$ such that $\Lambda_1(\{0\}) = \Lambda_2(\{0\}) = 0$ and $\Lambda_1(B) \geq \Lambda_2(B)$ for all measurable $B$. Then there exist $\mathcal{P}$-valued processes $(\Pi_1(t), t \geq 0)$ and $(\Pi_2(t), t \geq 0)$ such that $\Pi_1$ is a $\Lambda_1$-coalescent, $\Pi_2$ is a $\Lambda_2$-coalescent, and $N_1(t) \leq N_2(t)$ for all $t \geq 0$, where $N_i(t)$ is the number of blocks of $\Pi_i(t)$ for $i = 1, 2$.*

**Proof.** For $i = 1, 2$, let $L_i$ be the measure on $\{0,1\}^\infty$ such that $L_i(B) = \int_0^1 Q_x(B) x^{-2} \Lambda_i(dx)$ for all measurable $B$. Let $L_3(B) = L_1(B) - L_2(B) \geq 0$ for all measurable $B$. Let $\Psi_2$ and $\Psi_3$ be independent Poisson point processes with intensities $dt \times L_2(d\xi)$ and $dt \times L_3(d\xi)$, respectively. Let $\Psi_1$ be the Poisson point process consisting of all points in $\Psi_2$ and $\Psi_3$, which has intensity $dt \times L_1(d\xi)$. For $i = 1, 2$, let $(\Pi_i(t), t \geq 0)$ be the $\mathcal{P}$-valued coalescent process obtained from $\Psi_i$ as described above. Then $(\Pi_1(t), t \geq 0)$ is a $\Lambda_1$-coalescent and $(\Pi_2(t), t \geq 0)$ is a $\Lambda_2$-coalescent.

For $i = 1, 2$, let $N_{i,n}(t)$ be the number of blocks of $R_n \Pi_i(t)$. To show that $N_1(t) \leq N_2(t)$ for all $t \geq 0$, it suffices to show that $N_{1,n}(t) \leq N_{2,n}(t)$ for all positive integers $n$ and all $t \geq 0$. Each point of $\Psi_2$ is also a point of $\Psi_1$. Suppose $(t, \xi)$ is a point of $\Psi_2$, and let $A_t = \{i: \xi_i = 0 \text{ or } \xi_j = 0 \text{ for all } j < i\}$. Then $N_{1,n}(t)$ is the cardinality of $A_t \cap \{1, \ldots, N_{1,n}(t-)\}$ and $N_{2,n}(t)$ is the cardinality of $A_t \cap \{1, \ldots, N_{2,n}(t-)\}$. Therefore, if $N_{1,n}(t-) \leq N_{2,n}(t-)$, then $N_{1,n}(t) \leq N_{2,n}(t)$. Since $N_{1,n}(0) = N_{2,n}(0)$, the result follows from the construction and the fact that the restricted processes $(R_n \Pi_1(t), t \geq 0)$ and $(R_n \Pi_2(t), t \geq 0)$ have only finitely many jump times. □

## 4. Number of blocks

In this section, we prove Theorem 1.1. Throughout this section, we assume that $(M_t, t \geq 0)$ is the measure-valued process defined in Section 2 and that $R(t)$ and $R^{-1}(t)$ are as defined in Section 2. Since Lemma 2.2 gives the number of atoms of $M_t$, the key step involves analyzing the time change, which will make it possible to relate the number of blocks of the Beta$(2 - \alpha, \alpha)$-coalescent to the number of atoms of $M_t$.

**Lemma 4.1.** *Suppose $(Y(t), t \geq 0)$ is a Lévy process such that $Y(0) = 0$ and $E[e^{-\lambda Y(t)}] = e^{t\Psi(\lambda)}$, where $\Psi(\lambda) = \lambda^\alpha$ for some $\alpha \in (1, 2)$. There exists a constant $C$ such that for all $t > 0$ and $\varepsilon > 0$, we have*

$$P\left(\sup_{0 \leq s \leq t} |Y(s)| > \varepsilon\right) \leq Ct\varepsilon^{-\alpha}.$$

**Proof.** Scaling properties of Lévy process imply that for all $k > 0$, the processes $(Y(t), t \geq 0)$ and $(k^{-1/\alpha} Y(kt), t \geq 0)$ have the same law. By taking $k = 1/t$, we get

$$P\left(\sup_{0 \leq s \leq t} |Y(s)| > \varepsilon\right) = P\left(\sup_{0 \leq s \leq 1} |Y(s)| > t^{-1/\alpha} \varepsilon\right). \tag{15}$$

It is known for Lévy processes with the above scaling property (see Exercise 2 in Chapter VIII of [4]) that there exists a constant $C_1$ such that

$$P\Big(\sup_{0 \leq s \leq 1} |Y(s)| > x\Big) \sim C_1 x^{-\alpha}, \tag{16}$$

where $\sim$ means that the ratio of the two sides tends to 1 as $x \to \infty$. The lemma follows from (15) and (16). $\square$

**Lemma 4.2.** *There exists a constant $C$ such that for all $t > 0$ and $\varepsilon > 0$, we have*

$$P\Big(\frac{(1-\varepsilon)t}{\alpha(\alpha-1)\Gamma(\alpha)} \leq R^{-1}(t) \leq \frac{(1+\varepsilon)t}{\alpha(\alpha-1)\Gamma(\alpha)}\Big) \geq 1 - Ct\varepsilon^{-\alpha}.$$

**Proof.** Let $(Z(t), t \geq 0)$ be a CSBP with branching mechanism $\Psi(\lambda) = \lambda^\alpha$ such that $Z(0) = 1$. Since every CSBP can be obtained via a time change of a Lévy process with no negative jumps, as explained at the end of Section 2, we may assume that there is a Lévy process $(Y(t), t \geq 0)$ satisfying $Y(0) = 1$ and $E[e^{-\lambda Y(t)}] = e^{-\lambda + t\Psi(\lambda)}$ such that $Z(t) = Y(U(t))$ for all $t$, where $U(t) = \inf\{s : \int_0^s \tilde{Y}(u)^{-1} \, du > t\}$ and $(\tilde{Y}(t), t \geq 0)$ is the process $(Y(t), t \geq 0)$ stopped when it hits zero.

Assume $\varepsilon < 1/2$, and let $K = 4/[\alpha(\alpha-1)\Gamma(\alpha)]$. Assume that $|Y(s) - 1| \leq \varepsilon$ for all $s \in [0, Kt]$, which happens with probability at least $1 - Ct\varepsilon^{-\alpha}$ for some constant $C$ by Lemma 4.1. We then have $(1-\varepsilon)s \leq U(s) \leq (1+\varepsilon)s$ for all $s \in [0, Kt/2]$. Since

$$R(s) = \alpha(\alpha-1)\Gamma(\alpha) \int_0^s Y(U(r))^{1-\alpha} \, dr$$

for all $s$, it follows that

$$\alpha(\alpha-1)\Gamma(\alpha)(1+\varepsilon)^{1-\alpha} s \leq R(s) \leq \alpha(\alpha-1)\Gamma(\alpha)(1-\varepsilon)^{1-\alpha} s$$

for all $s \in [0, Kt/2]$. Therefore,

$$\frac{t(1-\varepsilon)^{\alpha-1}}{\alpha(\alpha-1)\Gamma(\alpha)} \leq R^{-1}(t) \leq \frac{t(1+\varepsilon)^{\alpha-1}}{\alpha(\alpha-1)\Gamma(\alpha)}.$$

Since $1 - \varepsilon \leq (1-\varepsilon)^{\alpha-1}$ and $(1+\varepsilon)^{\alpha-1} \leq 1 + \varepsilon$, the lemma follows. $\square$

**Lemma 4.3.** *Let $(\Pi(t), t \geq 0)$ be the Beta$(2-\alpha, \alpha)$ coalescent, where $1 < \alpha < 2$, and let $N(t)$ be the number of blocks of $\Pi(t)$. There exists a constant $C$ depending on $\varepsilon$ such that for all $t > 0$, we have*

$$P((1-\varepsilon)(\alpha\Gamma(\alpha))^{1/(\alpha-1)} \leq t^{1/(\alpha-1)} N(t) \leq (1+\varepsilon)(\alpha\Gamma(\alpha))^{1/(\alpha-1)}) \geq 1 - Ct.$$

**Proof.** Recall that $(M_t, t \geq 0)$ is the measure-valued process defined in Section 2. Let $D(s)$ be the number of atoms of $M_s$. By Lemma 2.2, the distribution of $D(s)$ is Poisson with mean $\theta_s = [(\alpha-1)s]^{-1/(\alpha-1)}$, so $E[D(s)] = \text{Var}(D(s)) = \theta_s$. By Chebyshev's Inequality, we have $P(|D(s) - \theta_s| > \delta\theta_s) \leq 1/(\delta^2 \theta_s)$. Therefore, if $\delta$ is small enough that $1 + \varepsilon \geq (1+\delta)(1-\delta)^{-1/(\alpha-1)}$ and $1 - \varepsilon \leq (1-\delta)(1+\delta)^{-1/(\alpha-1)}$, then

$$P\Big(D\Big(\frac{t(1-\delta)}{\alpha(\alpha-1)\Gamma(\alpha)}\Big) > (1+\varepsilon)(\alpha\Gamma(\alpha))^{1/(\alpha-1)} t^{-1/(\alpha-1)}\Big) \leq \frac{1}{\delta^2}\Big(\frac{1-\delta}{\alpha\Gamma(\alpha)}\Big)^{1/(\alpha-1)} t^{1/(\alpha-1)} \tag{17}$$

and likewise

$$P\Big(D\Big(\frac{t(1+\delta)}{\alpha(\alpha-1)\Gamma(\alpha)}\Big) < (1-\varepsilon)(\alpha\Gamma(\alpha))^{1/(\alpha-1)} t^{-1/(\alpha-1)}\Big) \leq \frac{1}{\delta^2}\Big(\frac{1+\delta}{\alpha\Gamma(\alpha)}\Big)^{1/(\alpha-1)} t^{1/(\alpha-1)}. \tag{18}$$

Since $N(t)$ has the same distribution as $D(R^{-1}(t))$ by Lemma 2.1 and $D(t)$ is a decreasing function of $t$, the result follows from (17), (18), and Lemma 4.2. □

We are now ready to prove Theorem 1.1 for the Beta$(2-\alpha,\alpha)$-coalescents.

**Proposition 4.4.** *Let $(\Pi(t), t \geq 0)$ be the Beta$(2-\alpha, \alpha)$-coalescent where $1 < \alpha < 2$, and let $N(t)$ be the number of blocks of $\Pi(t)$. Then*

$$\lim_{t \downarrow 0} t^{1/(\alpha-1)} N(t) = (\alpha \Gamma(\alpha))^{1/(\alpha-1)} \quad \text{a.s.}$$

**Proof.** Let $\varepsilon > 0$. Fix $t > 0$, and let $t_j = t(1-\varepsilon)^j$ for $j = 0, 1, 2, \ldots$. Let $B$ be the event that for all $j$, we have

$$(1-\varepsilon)(\alpha \Gamma(\alpha))^{1/(\alpha-1)} \leq t_j^{1/(\alpha-1)} N(t_j) \leq (1+\varepsilon)(\alpha \Gamma(\alpha))^{1/(\alpha-1)}. \tag{19}$$

By Lemma 4.3, there is a constant $C$ depending on $\varepsilon$ such that

$$P(B) \geq 1 - C \sum_{j=0}^{\infty} t(1-\varepsilon)^j = 1 - C\varepsilon^{-1} t. \tag{20}$$

Suppose $B$ occurs and $0 < s \leq t$. Then for some $j$, we have $t_{j+1} < s \leq t_j$, which implies that $N(t_j) \leq N(s) \leq N(t_{j+1})$, since the number of blocks is a decreasing function of time. From (19) and the definition of the $t_j$, we get

$$(1-\varepsilon)^{1+1/(\alpha-1)}(\alpha \Gamma(\alpha))^{1/(\alpha-1)} \leq s^{1/(\alpha-1)} N(s) \leq (1+\varepsilon)(1-\varepsilon)^{-1/(\alpha-1)}(\alpha \Gamma(\alpha))^{1/(\alpha-1)}.$$

Letting $t \downarrow 0$ and using (20), we get

$$\liminf_{t \downarrow 0} t^{1/(\alpha-1)} N(t) \geq (1-\varepsilon)^{1+1/(\alpha-1)}(\alpha \Gamma(\alpha))^{1/(\alpha-1)} \quad \text{a.s.}$$

and

$$\limsup_{t \downarrow 0} t^{1/(\alpha-1)} N(t) \leq (1+\varepsilon)(1-\varepsilon)^{-1/(\alpha-1)}(\alpha \Gamma(\alpha))^{1/(\alpha-1)} \quad \text{a.s.}$$

Letting $\varepsilon \downarrow 0$ completes the proof. □

**Remark 4.5.** *Note that one can not conclude Proposition 4.4 simply by combining the facts that that $D(t)$ is asymptotically equivalent to $\theta_t$ as $t \downarrow 0$ and that $R^{-1}(t)$ is asymptotically equivalent to $t/[\alpha(\alpha-1)\Gamma(\alpha)]$ as $t \downarrow 0$. The more involved argument in Proposition 4.4 is necessary because Lemma 2.2 only establishes an equality in distribution at individual times. Consequently, almost sure results about the CSBP as $t \downarrow 0$ do not immediately translate to the coalescent.*

**Proof of Theorem 1.1.** Let $\varepsilon > 0$, and then choose $\delta > 0$ such that $(A - \varepsilon)x^{1-\alpha}(1-x)^{\alpha-1} \leq f(x) \leq (A + \varepsilon)x^{1-\alpha}(1-x)^{\alpha-1}$ for all $x \in [0, \delta]$. Let $\Lambda_0$ be the finite measure on $[0,1]$ with density $f(x)\mathbf{1}_{\{x \leq \delta\}}$. Let $\Lambda_1$ be the finite measure on $[0,1]$ with density $(A-\varepsilon)x^{1-\alpha}(1-x)^{\alpha-1}\mathbf{1}_{\{x \leq \delta\}}$, and let $\Lambda_2$ be the finite measure on $[0,1]$ with density $(A+\varepsilon)x^{1-\alpha}(1-x)^{\alpha-1}\mathbf{1}_{\{x \leq \delta\}}$. For $i = 0, 1, 2$, let $(\overline{\Pi}_i(t), t \geq 0)$ be a $\Lambda_i$-coalescent, and let $N_i(t)$ denote the number of blocks in $\Pi_i(t)$.

Note that if $(\Pi(t), t \geq 0)$ is a $\Lambda'$-coalescent and $C$ is a constant, then $(\Pi(Ct), t \geq 0)$ is a $C\Lambda'$-coalescent. This fact, combined with Proposition 4.4 and Lemma 3.1, gives

$$\lim_{t \downarrow 0} t^{1/(\alpha-1)} N_1(t) = \left(\frac{\alpha}{(A-\varepsilon)\Gamma(2-\alpha)}\right)^{1/(\alpha-1)} \quad \text{a.s.}$$

and

$$\lim_{t \downarrow 0} t^{1/(\alpha-1)} N_2(t) = \left( \frac{\alpha}{(A+\varepsilon)\Gamma(2-\alpha)} \right)^{1/(\alpha-1)} \quad \text{a.s.}$$

By applying Lemma 3.2 to $\Lambda_0$ and $\Lambda_1$, we get

$$\limsup_{t \downarrow 0} t^{1/(\alpha-1)} N_0(t) \leq \left( \frac{\alpha}{(A-\varepsilon)\Gamma(2-\alpha)} \right)^{1/(\alpha-1)} \quad \text{a.s.} \tag{21}$$

Likewise, by applying Lemma 3.2 to $\Lambda_0$ and $\Lambda_2$, we get

$$\liminf_{t \downarrow 0} t^{1/(\alpha-1)} N_0(t) \geq \left( \frac{\alpha}{(A+\varepsilon)\Gamma(2-\alpha)} \right)^{1/(\alpha-1)} \quad \text{a.s.} \tag{22}$$

The conclusion of the theorem for the $\Lambda_0$-coalescent now follows by letting $\varepsilon \downarrow 0$ in (21) and (22). The conclusion for the original $\Lambda$-coalescent then follows from Lemma 3.1. $\square$

## 5. Block sizes

Our goal in this section is to prove Theorems 1.2 and 1.4 and Propositions 1.5 and 1.6, all of which pertain to the sizes of the blocks in the coalescent.

### 5.1. One-dimensional distributions

We first prove Theorem 1.2. Although Theorem 1.2 is stated for $\alpha \in (1, 2)$, the proof also works for $\alpha = 2$, so we get an alternative proof of the fact that for Kingman's coalescent, $\Theta(T_k)$ is uniformly distributed on $\Delta_k$.

**Proof of Theorem 1.2.** Let $\lambda_k = \sum_{j=2}^{k} \binom{k}{j} \lambda_{k,j}$ be the total rate of all mergers when the coalescent has $k$ blocks. Let $V_k = \{N(T_k) = k\}$ be the event that at some time the coalescent has exactly $k$ blocks. Conditional on $V_k$, the amount of time for which the coalescent has $k$ blocks has an exponential distribution with mean $\lambda_k^{-1}$, and this time is independent of $\Theta(T_k)$. Therefore, if $B$ is a measurable subset of $\Delta_k$, then

$$P(\Theta(T_k) \in B | V_k) = \frac{P(\{\Theta(T_k) \in B\} \cap V_k)}{P(V_k)} = \frac{\lambda_k}{P(V_k)} E\left[ \int_0^\infty \mathbf{1}_{\{N(t)=k, \Theta(t) \in B\}} \, dt \right]. \tag{23}$$

Let $(Z(t), t \geq 0)$ be a CSBP with branching mechanism $\Psi(\lambda) = \lambda^\alpha$, obtained from the measure-valued process $(M_t, t \geq 0)$ as in Section 2. Let $D(t)$ be the number of atoms of $M_t$, and let $J(t) = (J_1(t), \ldots, J_{D(t)}(t))$ be the sequence consisting of the sizes of the atoms of $M_t$, ranked in decreasing order. Let $J^*(t) = 0$ on $\{Z(t) = 0\}$, and let $J^*(t) = J(t)/Z(t)$ on $\{Z(t) > 0\}$, so the terms in the sequence $J^*(t)$ sum to one for all $t$ such that $Z(t) > 0$. By Lemma 2.1, the distribution of $(N(t), \Theta(t))$ is the same as the distribution of $(D(R^{-1}(t)), J^*(R^{-1}(t)))$. Combining this result with Fubini's theorem and then making the change of variables $s = R^{-1}(t)$, we have

$$E\left[ \int_0^\infty \mathbf{1}_{\{N(t)=k, \Theta(t) \in B\}} \, dt \right] = E\left[ \int_0^\infty \mathbf{1}_{\{D(R^{-1}(t))=k, J^*(R^{-1}(t)) \in B\}} \, dt \right]$$
$$= E\left[ \int_0^\infty \alpha(\alpha-1)\Gamma(\alpha) Z(s)^{1-\alpha} \mathbf{1}_{\{D(s)=k, J^*(s) \in B\}} \, ds \right]. \tag{24}$$

Recall from (14) that $\theta_s = [(\alpha-1)s]^{-1/(\alpha-1)}$ and therefore $\theta_s^{1-\alpha} = (\alpha-1)s$. Therefore, (23) and (24) imply

$$P(\Theta(T_k) \in B | V_k) = \frac{\lambda_k \alpha \Gamma(\alpha)}{P(V_k)} E\left[\int_0^\infty s^{-1}(\theta_s Z(s))^{1-\alpha} \mathbf{1}_{\{D(s)=k, J^*(s) \in B\}} \, ds\right]$$

$$= \frac{\lambda_k \alpha \Gamma(\alpha)}{P(V_k)} \int_0^\infty s^{-1} P(D(s)=k) E[(\theta_s Z(s))^{1-\alpha} \mathbf{1}_{\{J^*(s) \in B\}} | D(s)=k] \, ds.$$

Recall that $X_1, \ldots, X_k$ are obtained by picking $k$ i.i.d. random variables with distribution $\mu$, and then ranking them in decreasing order. Also, recall that $S_k = X_1 + \cdots + X_k$. By Lemma 2.2, the conditional distribution of $J(s)$ given $D(s) = k$ is the same as the distribution of $(\theta_s^{-1} X_1, \ldots, \theta_s^{-1} X_k)$. Because $Z(s) = J_1(s) + \cdots + J_{D(s)}(s)$, it follows that the joint distribution of $(\theta_s Z(s), J^*(s))$, conditional on $D(s) = k$, is the same as the joint distribution of $(S_k, (X_1/S_k, \ldots, X_k/S_k))$. Note that this distribution is the same for all $s$. Therefore, for $g = \mathbf{1}_B$, we have

$$P(\Theta(T_k) \in B | V_k) = \frac{\lambda_k \alpha \Gamma(\alpha)}{P(V_k)} \int_0^\infty s^{-1} P(D(s)=k) E\left[S_k^{1-\alpha} g\left(\frac{X_1}{S_k}, \ldots, \frac{X_k}{S_k}\right)\right] ds$$

$$= C E\left[S_k^{1-\alpha} g\left(\frac{X_1}{S_k}, \ldots, \frac{X_k}{S_k}\right)\right], \tag{25}$$

where $C = \lambda_k \alpha \Gamma(\alpha) P(V_k)^{-1} \int_0^\infty s^{-1} P(D(s)=k) \, ds$. By taking $B = \Delta_k$ so that both sides of (25) equal one, we get $C = E[S_k^{1-\alpha}]^{-1}$. This establishes (5) when $g$ is an indicator function. The result for arbitrary nonnegative measurable $g$ now follows from the linearity of expectation and the Monotone Convergence Theorem. □

*5.2. Block sizes at small times*

Our next goal is to prove Theorem 1.4. Our first lemma bounds the fluctuations of a continuous-state branching process for small times.

**Lemma 5.1.** *There exists a constant $C$ such that for all $a > 0$, $t > 0$, and $\varepsilon > 0$, if $(Z(t), t \geq 0)$ is a CSBP with stable branching mechanism $\Psi(\lambda) = \lambda^\alpha$ with $Z(0) = a$, then*

$$q(a, t, \varepsilon) = P\Big(\sup_{0 \leq s \leq t} |Z(s) - a| > \varepsilon\Big) \leq C(a+\varepsilon)t\varepsilon^{-\alpha}.$$

*In particular, for any constants $C_1$, $C_2$, and $C_3$, we have*

$$\lim_{t \to 0} q(C_1 t^{1/\alpha}, C_2 t, C_3 t^{1/\alpha}) = 0.$$

**Proof.** As in the proof of Lemma 4.2, we may assume that there is a Lévy process $(Y(t), t \geq 0)$ satisfying $Y_0 = a$ and $E[e^{-\lambda Y(t)}] = e^{-a\lambda + t\Psi(\lambda)}$ such that $Z(t) = Y(U(t))$ for all $t$, where $U(t) = \inf\{s: \int_0^s \tilde{Y}(u)^{-1} \, du > t\}$, and $(\tilde{Y}(t), t \geq 0)$ is the process $(Y(t), t \geq 0)$ stopped when it hits zero. If $|Y(s) - a| \leq \varepsilon$ for $0 \leq s \leq (a+\varepsilon)t$, then $U(s) \leq (a+\varepsilon)s$ for all $s \leq t$, and therefore, $|Z(s) - a| = |Y(U(s)) - a| \leq \varepsilon$ for all $s \leq t$. The result now follows from Lemma 4.1. □

**Lemma 5.2.** *Let $(\Pi(t), t \geq 0)$ be the Beta$(2-\alpha, \alpha)$ coalescent, where $1 < \alpha < 2$, and let $N(t, x)$ be the number of blocks of $\Pi(t)$ whose asymptotic frequency is at most $x$. Let $\gamma = (\alpha \Gamma(\alpha))^{1/(\alpha-1)}$. There exists a constant $C$ depending on $\varepsilon$ and $x$ such that for all $t > 0$, we have*

$$P((1-\varepsilon)\gamma F((1-\varepsilon)\gamma x) \leq t^{1/(\alpha-1)} N(t, t^{1/(\alpha-1)} x) \leq (1+\varepsilon)\gamma F((1+\varepsilon)\gamma x)) \geq 1 - Ct^{1/2}.$$

**Proof.** Fix $t > 0$ and $\varepsilon > 0$. Define $t_- = (1 - t^{1/2\alpha})t/[\alpha(\alpha-1)\Gamma(\alpha)]$ and $t_+ = (1 + t^{1/2\alpha})t/[\alpha(\alpha-1)\Gamma(\alpha)]$. Let $B_{1,t}$ be the event that $t_- \leq R^{-1}(t) \leq t_+$. Let $B_{2,t}$ be the event that $|Z(s) - 1| \leq t^{1/2\alpha}$ for all $s \leq t_+$. Let $A_3 = \{a\colon 0 < M_s(\{a\}) \leq (1 - t^{1/2\alpha})t^{1/(\alpha-1)}x \text{ for all } s \in [t_-, t_+]\}$, and let $A_4 = \{a\colon 0 < M_s(\{a\}) \leq (1 + t^{1/2\alpha})t^{1/(\alpha-1)}x \text{ for some } s \in [t_-, t_+]\}$. Then, letting $\#S$ denote the cardinality of the set $S$, we define $B_{3,t}$ to be the event that $\#A_3 \geq (1-\varepsilon)t^{-1/(\alpha-1)}\gamma F((1-\varepsilon)\gamma x)$ and $B_{4,t}$ to be the event that $\#A_4 \leq (1+\varepsilon)t^{-1/(\alpha-1)}\gamma F((1+\varepsilon)\gamma x)$.

By Lemmas 2.1 and 2.2, the distribution of $N(t, t^{1/(\alpha-1)}x)$ is the same as the distribution of the number of terms of the sequence $J(R^{-1}(t))/Z(R^{-1}(t))$ that are in $(0, t^{1/(\alpha-1)}x]$. Furthermore, note that if $B_{i,t}$ occurs for $i = 1, \ldots, 4$, then the number of terms of $J(R^{-1}(t))/Z(R^{-1}(t))$ in $(0, t^{1/(\alpha-1)}x]$ is at least $(1 - \varepsilon)t^{-1/(\alpha-1)}\gamma F((1-\varepsilon)\gamma x)$ and at most $(1+\varepsilon)t^{-1/(\alpha-1)}\gamma F((1+\varepsilon)\gamma x)$. Therefore, to prove the lemma, it suffices to show that there exists a constant $C$ such that $P(B_{1,t} \cap B_{2,t} \cap B_{3,t} \cap B_{4,t}) \geq 1 - Ct^{1/2}$. We get $P(B_{1,t}) \geq 1 - Ct^{1/2}$ from Lemma 4.2 and $P(B_{2,t}) \geq 1 - Ct^{1/2}$ from Lemma 5.1. It remains to consider $B_{3,t}$ and $B_{4,t}$.

We first bound $P(B_{3,t})$. By Lemma 2.2, if we let $\theta_t = [(\alpha-1)t]^{-1/(\alpha-1)}$, then for all $t > 0$ and $x > 0$, the number of atoms of $M_t$ with size at most $\theta_t^{-1}x$ has the Poisson distribution with mean $\theta_t F(x)$. Therefore, the number of atoms of $M_{t_-}$ with size at most $(1-\varepsilon)t^{1/(\alpha-1)}x$ has the Poisson distribution with mean

$$\theta_{t_-} F(\theta_{t_-}(1-\varepsilon)t^{1/(\alpha-1)}x) = \left(\frac{(1-t^{1/2\alpha})t}{\alpha\Gamma(\alpha)}\right)^{-1/(\alpha-1)} F((1-t^{1/2\alpha})^{-1/(\alpha-1)}(1-\varepsilon)\gamma x)$$

$$\geq t^{-1/(\alpha-1)}\gamma F((1-\varepsilon)\gamma x).$$

For sufficiently small $t$, we have $(1-t^{1/2\alpha})t^{1/(\alpha-1)}x - (1-\varepsilon)t^{1/(\alpha-1)}x \geq (\varepsilon/2)t^{1/(\alpha-1)}x$. For such $t$, the Markov property implies that conditional on $0 < M_{t_-}(\{a\}) \leq (1-\varepsilon)t^{1/(\alpha-1)}x$, the probability that $0 < M_s(\{a\}) < (1-t^{1/2\alpha})t^{1/(\alpha-1)}x$ for all $s \in [t_-, t_+]$ is at least

$$1 - q\left((1-\varepsilon)t^{1/(\alpha-1)}x, t_+ - t_-, \left(\frac{\varepsilon}{2}\right)t^{1/(\alpha-1)}x\right),$$

where $q$ is the function defined in Lemma 5.1. By Lemma 5.1, there is a constant $C$ such that

$$q\left((1-\varepsilon)t^{1/(\alpha-1)}x, t_+ - t_-, \left(\frac{\varepsilon}{2}\right)t^{1/(\alpha-1)}x\right) \leq Ct^{1/2\alpha}\varepsilon^{-\alpha},$$

which for sufficiently small $t$ is at most $\varepsilon/2$. Thus, for sufficiently small $t$, the cardinality of $A_3$ has the Poisson distribution with mean at least

$$\left(1 - \frac{\varepsilon}{2}\right)t^{-1/(\alpha-1)}\gamma F((1-\varepsilon)\gamma x).$$

Chebyshev's Inequality now implies that there is a constant $C$ such that for sufficiently small $t$, we have $P(B_{3,t}) \geq 1 - C\varepsilon^{-2}t^{1/(\alpha-1)}$.

We now need to bound $P(B_{4,t})$. The number of atoms of $M_{t_+}$ with size at most $(1+\varepsilon)t^{1/(\alpha-1)}x$ has the Poisson distribution with mean

$$\theta_{t_+} F(\theta_{t_+}(1+\varepsilon)t^{1/(\alpha-1)}x) = \left(\frac{(1+t^{1/2\alpha})t}{\alpha\Gamma(\alpha)}\right)^{-1/(\alpha-1)} F((1+t^{1/2\alpha})^{-1/(\alpha-1)}(1+\varepsilon)\gamma x)$$

$$\leq t^{-1/(\alpha-1)}\gamma F((1+\varepsilon)\gamma x).$$

For sufficiently small $t$, we have $(1+\varepsilon)t^{1/(\alpha-1)}x - (1+t^{1/2\alpha})t^{1/(\alpha-1)}x \geq (\varepsilon/2)t^{1/(\alpha-1)}x$. For every value of $a$ such that $0 < M_s(\{a\}) \leq (1+t^{1/2\alpha})t^{1/(\alpha-1)}x$ for some $s \in [t_-, t_+]$, we can apply the strong Markov property at the time $\inf\{s \geq t_-\colon M_s(\{a\}) \leq (1+t^{1/2\alpha})t^{1/(\alpha-1)}x\}$ to see that conditional on $0 < M_s(\{a\}) \leq (1+t^{1/2\alpha})t^{1/(\alpha-1)}x$ for some $s \in [t_-, t_+]$, the probability that $0 < M_{t_+}(\{a\}) \leq (1+\varepsilon)t^{1/(\alpha-1)}x$ is at least

$$1 - q\left((1+t^{1/2\alpha})t^{1/(\alpha-1)}x, t_+ - t_-, \left(\frac{\varepsilon}{2}\right)t^{1/(\alpha-1)}x\right) \geq 1 - Ct^{1/2\alpha}\varepsilon^{-\alpha} \geq \frac{1}{1+\varepsilon/2}$$

for sufficiently small $t$. It follows that for sufficiently small $t$, the cardinality of $A_{t,4}$ has a Poisson distribution with mean at most

$$\left(1 + \frac{\varepsilon}{2}\right)t^{-1/(\alpha-1)}\gamma F((1+\varepsilon)\gamma x).$$

The desired lower bound on $P(B_{t,4})$ now follows as before from Chebyshev's Inequality. □

**Proof of Theorem 1.4.** The proof is now similar to the proof of Proposition 4.4. Fix $x \geq 0$. We will first show that

$$\lim_{t \downarrow 0} t^{1/(\alpha-1)} N(t, t^{1/(\alpha-1)}x) = (\alpha\Gamma(\alpha))^{1/(\alpha-1)} F((\alpha\Gamma(\alpha))^{1/(\alpha-1)}x) \quad \text{a.s.} \tag{26}$$

Let $\varepsilon > 0$, and let $t > 0$. Let $t_j = t(1-\varepsilon)^j$ for $j = 0, 1, 2, \ldots$. Let $B$ be the event that for all $j$, we have

$$(1-\varepsilon)\gamma F((1-\varepsilon)\gamma x) \leq t_j^{1/(\alpha-1)} N(t_j, t_j^{1/(\alpha-1)}x) \leq (1+\varepsilon)\gamma F((1+\varepsilon)\gamma x), \tag{27}$$

where $\gamma = (\alpha\Gamma(\alpha))^{1/(\alpha-1)}$ as in Lemma 5.2. By Lemma 5.2, we have $P(B) \geq 1 - Ct^{1/2}$ for some constant $C$ which depends on $\varepsilon$.

The number of blocks in the coalescent with asymptotic frequency at most $x$ can only decrease as a result of mergers. Therefore, for each fixed $x$, $N(t,x)$ is a decreasing function of $t$, so if $t_{j+1} < s \leq t_j$, then $N(t_j, s^{1/(\alpha-1)}x) \leq N(s, s^{1/(\alpha-1)}x) \leq N(t_{j+1}, s^{1/(\alpha-1)}x)$. It follows that if $B$ occurs, then

$$N(t_j, t_j^{1/(\alpha-1)}(1-\varepsilon)^{1/(\alpha-1)}x) \leq N(s, s^{1/(\alpha-1)}x) \leq N(t_{j+1}, t_{j+1}^{1/(\alpha-1)}(1-\varepsilon)^{-1/(\alpha-1)}x).$$

Let $B'$ be the event that $B$ occurs and also that (27) holds for all $j \geq 0$ with $(1-\varepsilon)^{1/(\alpha-1)}x$ and $(1-\varepsilon)^{-1/(\alpha-1)}x$ in place of $x$. Note that $P(B') \geq 1 - C't^{1/2}$, where $C'$ is a constant which depends on $\varepsilon$. Now if $B'$ occurs, then for all $s \leq t$, we have

$$(1-\varepsilon)^{1+1/(\alpha-1)}\gamma F((1-\varepsilon)^{1+1/(\alpha-1)}\gamma x) \leq s^{1/(\alpha-1)}N(s, s^{1/(\alpha-1)}x)$$
$$\leq (1+\varepsilon)(1-\varepsilon)^{-1/(\alpha-1)}\gamma F((1+\varepsilon)(1-\varepsilon)^{-1/(\alpha-1)}\gamma x).$$

To obtain (26) by letting $t \downarrow 0$ and then $\varepsilon \downarrow 0$ as in the proof of Proposition 4.4, it remains only to show that $F$ is continuous or, equivalently, that $\mu$ has no atoms. This was proved in [8] for a measure that can be obtained by a rescaling of $\mu$, and we use the same argument here. Suppose $\mu(\{b\}) > 0$. By Lemma 2.2, we have $P(D(t) = 1) > 0$ for all $t > 0$, and therefore $[(\alpha-1)t]^{1/(\alpha-1)}b$ is an atom of the distribution of $Z(t)$. It then follows by applying the Markov property at time $1-t$ that $[(\alpha-1)t]^{1/(\alpha-1)}b$ is an atom of the distribution of $Z(1)$ for all $t \in (0,1]$, which is a contradiction.

It remains to establish that the convergence in (26) is uniform in $x$. Let $\varepsilon > 0$, and choose $N > 2/\varepsilon$. Since $F$ is continuous, we can choose $x_1, \ldots, x_{N-1}$ such that $F((\alpha\Gamma(\alpha))^{1/(\alpha-1)}x_i) = i/N$ for all $i$. Also set $x_0 = 0$ and $x_N = \infty$. By (26), almost surely for sufficiently small $t$ we have $|(\alpha\Gamma(\alpha))^{-1/(\alpha-1)}t^{1/(\alpha-1)}N(t, t^{1/(\alpha-1)}x_i) - i/N| < \varepsilon/2$ for $i = 0, \ldots, N$. For such $t$, we have

$$\sup_{x \geq 0}|(\alpha\Gamma(\alpha))^{-1/(\alpha-1)}t^{1/(\alpha-1)}N(t, t^{1/(\alpha-1)}x) - F((\alpha\Gamma(\alpha))^{1/(\alpha-1)}x)| < \varepsilon, \tag{28}$$

and the lemma follows. □

### 5.3. Size of the block containing 1

We now work towards proving Proposition 1.5, which concerns the distribution of the size of the block containing 1 or, equivalently, the distribution of the size of a size-biased pick from the blocks of the coalescent. We will deduce Proposition 1.5 from Theorem 1.4. We first review some facts about size-biased distributions. If $X = (X_1, X_2, \ldots)$ is a sequence of nonnegative random variables whose sum is 1, then a size-biased pick from the sequence $X$ is a random variable $X_N$ such that $P(N = n|X) = X_n$. If $(\Theta(t), t \geq 0)$ is a ranked $\Lambda$-coalescent with proper frequencies and $K(t)$ is the size of the block containing 1 at time $t$, then $K(t)$ is a size-biased pick from the sequence $\Theta(t)$.

If $X$ is a nonnegative random variable with finite mean, then the size-biased distribution of $X$ is the distribution of the random variable $\hat{X}$, where

$$E[f(\hat{X})] = \frac{E[Xf(X)]}{E[X]} \tag{29}$$

for all nonnegative measurable functions $f$. The next lemma records two facts about size-biased distributions that we will use in the proof of Proposition 1.5.

**Lemma 5.3.** *Suppose $X$ is a nonnegative random variable with mean 1. Let $\hat{X}$ be a random variable having the size-biased distribution of $X$. Then, for all $y \geq 0$,*

$$P(\hat{X} \leq y) = \int_0^y (P(X \leq y) - P(X \leq x)) \, dx. \tag{30}$$

*Let $\phi(\lambda) = E[e^{-\lambda X}]$. Then for all $\lambda > 0$,*

$$E[e^{-\lambda \hat{X}}] = -\phi'(\lambda). \tag{31}$$

**Proof.** Let $\mu$ denote the distribution of $X$. By (29) and Fubini's theorem,

$$P(\hat{X} \leq y) = E[X \mathbf{1}_{\{X \leq y\}}] = \int_0^y z \mu(dz) = \int_0^y \int_0^z dx \, \mu(dz) = \int_0^y \int_x^y \mu(dz) \, dx,$$

which leads to (30). To prove (31), note that (29) gives $E[e^{-\lambda \hat{X}}] = E[Xe^{-\lambda X}]$. It is also easily verified that $-\phi'(\lambda) = E[Xe^{-\lambda X}]$, where we can use the dominated convergence theorem to interchange differentiation and expectation because $E[X] < \infty$. $\square$

**Proof of Proposition 1.5.** For $t > 0$ and $\varepsilon > 0$, let $A_{t,\varepsilon}$ be the event that (28) occurs. Let $(\Theta(t), t \geq 0)$ be the ranked coalescent process associated with $(\Pi(t), t \geq 0)$, and write $\Theta(t) = (\Theta_1(t), \Theta_2(t), \ldots)$. Since $K(t)$ is a size-biased pick from $\Theta(t)$, we have, for all $y \geq 0$,

$$P(K(t) \leq y | A_{t,\varepsilon}) = E\left[\sum_{i=1}^\infty \Theta_i(t) \mathbf{1}_{\{\Theta_i(t) \leq y\}} \Big| A_{t,\varepsilon}\right] = E\left[\sum_{i=1}^\infty \int_0^y \mathbf{1}_{\{x < \Theta_i(t) \leq y\}} \, dx \Big| A_{t,\varepsilon}\right]$$

$$= E\left[\int_0^y (N(t,y) - N(t,x)) \, dx \Big| A_{t,\varepsilon}\right].$$

Therefore, letting $\gamma = (\alpha \Gamma(\alpha))^{1/(\alpha-1)}$, we get

$$P(\gamma t^{-1/(\alpha-1)} K(t) \leq y | A_{t,\varepsilon}) = E\left[\int_0^{\gamma^{-1} t^{1/(\alpha-1)} y} (N(t, \gamma^{-1} t^{1/(\alpha-1)} y) - N(t,x)) \, dx \Big| A_{t,\varepsilon}\right]$$

$$= E\left[\int_0^y \gamma^{-1} t^{1/(\alpha-1)} (N(t, \gamma^{-1} t^{1/(\alpha-1)} y) - N(t, \gamma^{-1} t^{1/(\alpha-1)} z)) \, dz \Big| A_{t,\varepsilon}\right].$$

Let $Y$ be a random variable with the distribution $\mu$, and let $\hat{Y}$ have the size-biased distribution of $Y$. By applying the definition of $A_{t,\varepsilon}$ and then Lemma 5.3, we get that there is a number $\theta$ such that $-2y\varepsilon \leq \theta \leq 2y\varepsilon$ and

$$P(\gamma t^{-1/(\alpha-1)} K(t) \leq y | A_{t,\varepsilon}) = \theta + \int_0^y (F(y) - F(z))\,dz = \theta + P(\hat{Y} \leq y).$$

Now, fix $\delta > 0$. Choose $\varepsilon$ small enough that $2y\varepsilon < \delta/2$, and then choose $t$ small enough that $P(A_{t,\varepsilon}) < \delta/2$, which is possible by Theorem 1.4. Then

$$|P(\gamma t^{-1/(\alpha-1)} K(t) \leq y) - P(\hat{Y} \leq y)| < \delta.$$

It follows that as $t \downarrow 0$, we have $\gamma t^{-1/(\alpha-1)} K(t) \to_d \hat{Y}$. The proposition now follows from Lemma 5.3, as the formula for the Laplace transform of $X$ in Proposition 1.5 comes from differentiating the right-hand side of (4). □

*5.4. The largest block*

Next we prove Proposition 1.6. This will require understanding the tail of the distribution $\mu$. A key tool will be the following Tauberian theorem, which comes from Theorem 8.1.6 of [10].

**Lemma 5.4.** *Let $X$ be a nonnegative random variable. For nonnegative integers $n$, let $\mu_n = E[X^n]$. For $\lambda > 0$, let $\phi(\lambda) = E[e^{-\lambda X}]$. If $\mu_n < \infty$, let $g_n(\lambda) = \mu_n - (-1)^n \phi^{(n)}(\lambda)$, where $\phi^{(n)}$ denotes the $n$th derivative of $\phi$. Suppose $L$ is a function that is slowly varying at infinity. If $\gamma = n + \beta$, where $0 < \beta < 1$, then the following are equivalent:*

$$g_n(\lambda) \sim \frac{\Gamma(\gamma+1)}{\Gamma(\beta+1)} \lambda^\beta L\left(\frac{1}{\lambda}\right) \quad \text{as } \lambda \downarrow 0, \tag{32}$$

$$P(X > x) \sim \frac{(-1)^n x^{-\gamma} L(x)}{\Gamma(1-\gamma)} \quad \text{as } x \to \infty,$$

*where $\sim$ means that the ratio of the two sides tends to one.*

This leads to the following result concerning the largest atom of $M_t$.

**Lemma 5.5.** *Let $(M_t, t \geq 0)$ be the measure-valued process defined in Section 2. Let $J_1(t)$ be the size of the largest atom of $M_t$. Then, for all $x > 0$, we have*

$$\lim_{t \downarrow 0} P(J_1(t) \leq t^{1/\alpha} x) = e^{-(\alpha-1)x^{-\alpha}/\Gamma(2-\alpha)}.$$

**Proof.** Let $X$ be a random variable with distribution $\mu$, where the Laplace transform of $\mu$ is given by (4). Since $E[X] = 1$, we can apply Lemma 5.4 with $n = 1$. Defining $\phi$ and $g_1$ as in Lemma 5.4, we get $\phi'(\lambda) = -(1+\lambda^{\alpha-1})^{-\alpha/(\alpha-1)}$ and therefore

$$g_1(\lambda) = 1 - (1+\lambda^{\alpha-1})^{-\alpha/(\alpha-1)} \sim \frac{\alpha}{\alpha-1} \lambda^{\alpha-1},$$

where $\sim$ means that the ratio of the two sides tends to one as $\lambda \downarrow 0$. It follows that (32) holds if $L(x) = 1/(\alpha-1)$ for all $x$. Therefore, by Lemma 5.4, as $x \to \infty$ we have

$$P(X > x) \sim \frac{(-1)x^{-\alpha}}{(\alpha-1)\Gamma(1-\alpha)} = \frac{x^{-\alpha}}{\Gamma(2-\alpha)}. \tag{33}$$

By Lemma 2.2, the number of atoms of $M_t$ is Poisson with mean $\theta_t = [(\alpha - 1)t]^{-1/(\alpha-1)}$, and the sizes of the atoms have the distribution of $\theta_t^{-1} X$, where $X$ has distribution $\mu$. Therefore, for any $y > 0$, the distribution of the number of atoms of size at least $y$ is Poisson with mean $\theta_t P(X > \theta_t y)$. By (33),

$$\lim_{t \downarrow 0} \theta_t P(X > \theta_t t^{1/\alpha} x) = \frac{\theta_t (\theta_t t^{1/\alpha} x)^{-\alpha}}{\Gamma(2-\alpha)} = \frac{x^{-\alpha}(\alpha-1)}{\Gamma(2-\alpha)},$$

which implies the result. □

**Proof of Proposition 1.6.** Let $(Z(t), t \geq 0)$ be a CSBP obtained from the measure-valued process $(M_t, t \geq 0)$, as defined in Section 2. Let $J_1(t)$ be the size of the largest atom of $M_t$. The distribution of $W(t)$ is the same as the distribution of $J_1(R^{-1}(t))/Z(R^{-1}(t))$ by Lemma 2.1. By Lemma 4.2 and the right-continuity of $Z$, we have $Z(R^{-1}(t)) \to 1$ a.s. as $t \to \infty$. Therefore, it suffices to show that for all $x > 0$, if $\gamma = (\alpha \Gamma(\alpha) \Gamma(2-\alpha))^{-1/\alpha}$, then

$$\lim_{t \downarrow 0} P(J_1(R^{-1}(t)) \leq \gamma t^{1/\alpha} x) = e^{-x^{-\alpha}}. \tag{34}$$

We follow a strategy similar to that used in the proof of Lemma 5.2. Let $\varepsilon > 0$ and $t > 0$. Let $t_- = (1-\varepsilon)t/[\alpha(\alpha-1)\Gamma(\alpha)]$ and let $t_+ = (1+\varepsilon)t/[\alpha(\alpha-1)\Gamma(\alpha)]$. Let $B_{1,t}$ be the event that $t_- \leq R^{-1}(t) \leq t_+$. Let $B_{2,t}$ be the event that for some $a \in [0,1]$, we have $M_s(\{a\}) > \gamma t^{1/\alpha} x$ for all $s \in [t_-, t_+]$. Let $B_{3,t}$ be the event that $J_1(s) < \gamma t^{1/\alpha} x$ for all $s \in [t_-, t_+]$. Note that on the event $B_{1,t} \cap B_{2,t}$, we have $J_1(R^{-1}(t)) > \gamma t^{1/\alpha} x$, while on the event $B_{1,t} \cap B_{3,t}$, we have $J_1(R^{-1}(t)) < \gamma t^{1/\alpha} x$. Also, note that $\lim_{t \downarrow 0} P(B_{1,t}) = 1$ by Lemma 4.2.

Recall the definition of $q(a, t, \varepsilon)$ from Lemma 5.1. Using the Markov property for CSBPs at time $t_-$, we get

$$P(B_{2,t}) \geq P(J_1(t_-) > (1+\varepsilon)\gamma t^{1/\alpha} x)(1 - q((1+\varepsilon)\gamma t^{1/\alpha} x, t_+ - t_-, \varepsilon \gamma t^{1/\alpha} x)). \tag{35}$$

By Lemma 5.5 and the definition of $t_-$,

$$\lim_{t \downarrow 0} P(J_1(t_-) > (1+\varepsilon)\gamma t^{1/\alpha} x) = 1 - e^{-x^{-\alpha}(1-\varepsilon)(1+\varepsilon)^{-\alpha}}. \tag{36}$$

By Lemma 5.1,

$$\lim_{t \downarrow 0} q((1+\varepsilon)\gamma t^{1/\alpha} x, t_+ - t_-, \varepsilon \gamma t^{1/\alpha} x) = 0. \tag{37}$$

Combining (35), (36) and (37), we get

$$\limsup_{t \downarrow 0} P(J_1(R^{-1}(t)) \leq \gamma t^{1/\alpha} x) \leq e^{-x^{-\alpha}(1-\varepsilon)(1+\varepsilon)^{-\alpha}}. \tag{38}$$

To estimate $P(B_{3,t})$, we first define $B_{4,t}$ to be the event that $J_1(t_+) \leq (1-\varepsilon)\gamma t^{1/\alpha} x$. By applying the strong Markov property at the stopping time $T = \inf\{s \geq t_- : J_1(s) \geq \gamma t^{1/\alpha} x\}$, we see that

$$P(B_{4,t}|B_{3,t}^c) \leq q((1-\varepsilon)\gamma t^{1/\alpha} x, t_+ - t_-, \varepsilon \gamma t^{1/\alpha} x). \tag{39}$$

Therefore,

$$P(B_{3,t}) \geq P(B_{4,t} \cap B_{3,t}) = P(B_{4,t}) - P(B_{4,t} \cap B_{3,t}^c) \geq P(B_{4,t}) - P(B_{4,t}|B_{3,t}^c),$$

which by (39) and Lemma 5.1 converges to $P(B_{4,t})$ as $t \downarrow 0$. Combining this result with Lemma 5.5, we get

$$\liminf_{t \downarrow 0} P(J_1(R^{-1}(t)) \leq \gamma t^{1/\alpha} x) \geq \liminf_{t \downarrow 0} P(B_{4,t}) = e^{-x^{-\alpha}(1+\varepsilon)(1-\varepsilon)^{-\alpha}}. \tag{40}$$

Combining (38) and (40) and letting $\varepsilon \downarrow 0$ gives (34). □

## 6. Hausdorff and packing dimensions

In this section, we prove Theorem 1.7. Because it is known that $\dim_H(S) \leq \dim_P(S)$, it suffices to show that $\dim_P(S) \leq 1/(\alpha - 1)$ a.s. and that $\dim_H(S) \geq 1/(\alpha - 1)$ a.s. We begin with the upper bound on the packing dimension, which is easier and follows from the same argument given by Evans [20] for Kingman's coalescent.

**Proposition 6.1.** *We have $\dim_P(S) \leq 1/(\alpha - 1)$ a.s.*

**Proof.** By (8), it suffices to show that for all $\beta > 1/(\alpha - 1)$, we have $p_\beta(S) < \infty$, where $p_\beta(S)$ is as defined in (7). Since $S$ itself is a cover of $S$, it suffices to show that $P_\beta(S) < \infty$ for all $\beta > 1/(\alpha - 1)$, where $P_\beta(S)$ is as defined in (6). For this it suffices to show that for all $\beta > 1/(\alpha - 1)$, almost surely we can find $\delta > 0$ and $K < \infty$ such that if $\{V_i\}_{i=1}^\infty$ is a $\delta$-packing of $S$, then $\sum_{i=1}^\infty |V_i|^\beta \leq K < \infty$.

Let $N(t)$ be the number of blocks of $\Pi(t)$, and let $\beta > 1/(\alpha - 1)$. By Theorem 1.1, almost surely there exists a $\delta > 0$ such that for all $t \in (0, \delta)$, we have $N(t) \leq Ct^{-1/(\alpha-1)}$, where $C = 2[\alpha/(A\Gamma(2-\alpha))]^{1/(\alpha-1)}$. Let $\{V_i\}_{i=1}^\infty$ be a $\delta$-packing of $S$, and let $D_k = \#\{i: 2^{-(k+1)}\delta < |V_i| \leq 2^{-k}\delta\}$. For all $x \in S$, let $B(x, r)$ denote the closed ball of radius $r$ centered at $x$. Note that this ball also has diameter $r$ because $(S, d)$ is an ultrametric space. By construction, for all $t > 0$, the set $S$ is a union of $N(t)$ balls of the form $B(x, t)$, each containing the integers in one of the blocks of $\Pi(t)$. If $s > t$, then any open ball of radius $s$ must have its center in one of the $N(t)$ closed balls of radius $t$, and therefore must contain one of the $N(t)$ balls of radius $t$. It follows that $D_k \leq N(2^{-(k+1)}\delta)$. Thus,

$$\sum_{i=1}^\infty |V_i|^\beta \leq \sum_{k=0}^\infty D_k(2^{-k}\delta)^\beta \leq \sum_{k=0}^\infty N(2^{-(k+1)}\delta)(2^{-k}\delta)^\beta \leq C\delta^{-1/(\alpha-1)+\beta} \sum_{k=0}^\infty 2^{(k+1)/(\alpha-1)-k\beta} \leq K < \infty,$$

where $K < \infty$ because $\beta > 1/(\alpha - 1)$. $\square$

It remains to prove a lower bound on the Hausdorff dimension. Our argument is motivated by the proof of Theorem 5.5 in [17]. We will use the following lemma, which is essentially Proposition 4.9a in [21]. Proposition 4.9a in [21] is stated for Euclidean space, but the same proof works in general metric spaces.

**Lemma 6.2.** *Let $\gamma$ be a probability measure on $S$. If*

$$\limsup_{r \downarrow 0} \frac{\gamma(B(x, r))}{r^\beta} < C \quad \text{for } \gamma\text{-almost all } x \in S,$$

*then $\dim_H(S) \geq \beta$.*

To apply Lemma 6.2, it will be necessary to construct a suitable probability measure $\gamma$ on $S$. We will use the same approach used by Evans for Kingman's coalescent in Section 5 of [20]. If $B(x, t)$ is a ball in $S$, then $B(x, t)$ contains the integers in one of the $N(t)$ blocks of $\Pi(t)$, and all of these integers are centers of the ball. Consequently, every ball in $S$ can be written as $B(n, t)$, where $n \in \mathbb{N}$, and since the coalescent process $(\Pi(t), t \geq 0)$ has jumps only at a countable set of times, only countably many of these balls are distinct. Given $n \in \mathbb{N}$ and $t > 0$, define $\gamma(B(n, t))$ to be the asymptotic frequency of the block of $\Pi(t)$ containing $n$. If $s > t$, then the block $B$ of $\Pi(s)$ containing $n$ is a union of finitely many blocks $B_1, \ldots, B_k$ of $\Pi(t)$, and if $n_i \in B_i$ for $i = 1, \ldots, k$, then $B(n, s) = \bigcup_{i=1}^k B(n_i, t)$. Since the asymptotic frequency of $B$ is the sum of the asymptotic frequencies of $B_1, \ldots, B_k$, the function $\gamma$ can easily be extended to a finitely additive set function on the collection $\mathcal{B}$ consisting of the finite unions of balls. One can easily check that the complement of a finite union of balls in $S$ is also a finite union of balls, so $\mathcal{B}$ is an algebra. Every open subset of $S$ can be written as a union of balls, and therefore as a countable union of balls centered at one of the integers, so the $\sigma$-algebra generated by $\mathcal{B}$ is the Borel $\sigma$-algebra. Since $S$ is complete and almost surely can be covered by finitely many balls of radius $t$ for any $t > 0$, we have that $S$ is compact almost surely. As noted in [20], it then follows from Theorems 3.1.1 and 3.1.4 in [15] that $\gamma$ can be extended to a probability measure on $S$.

**Lemma 6.3.** *Let $K(t)$ be the asymptotic frequency of the block of $\Pi(t)$ containing $1$. Let $X$ be a point of $S$ chosen at random with distribution $\gamma$. Then the processes $(K(t), t \geq 0)$ and $(\gamma(B(X,t)), t \geq 0)$ have the same law.*

**Proof.** Clearly for all positive integers $i$ and $n$, the indicator of the event that $i$ and $n$ are in the same block of $\Pi(t)$ is measurable. Therefore, $\gamma(B(n,t))$ is measurable for all positive integers $n$, as it can be expressed as a limit of an average of these indicators. It follows that $\gamma(B(X,t))$ is also measurable. The processes $(K(t), t \geq 0)$ and $(\gamma(B(X,t)), t \geq 0)$ are both nondecreasing, so it suffices to show that they have the same finite-dimensional distributions. Since $K(0) = \gamma(B(X,0)) = 0$ a.s., it suffices to show that, for all $t > 0$, the processes $(K(s), s \geq t)$ and $(\gamma(B(X,s)), s \geq t)$ have the same law.

Let $(\Pi'(t), t \geq 0)$ be the restriction of $(\Pi(t), t \geq 0)$ to $\mathbb{N}' = \{2, 3, \ldots\}$, meaning that if $i, j \geq 2$, then $i$ and $j$ are in the same block of $\Pi'(t)$ if and only if they are in the same block of $\Pi$. Let $d'$ be the restriction of the metric $d$ to $\mathbb{N}'$, and let $(S', d')$ be the completion of $(\mathbb{N}', d')$. Since the $\Lambda$-coalescent has proper frequencies, the integer $1$ is not a singleton in $\Pi(t)$ for any $t$. Therefore, if $(t_k)_{k=1}^\infty$ is a sequence of positive numbers converging to zero, for each $k$ there is an integer $n_k$ such that $1$ and $n_k$ are in the same block of $\Pi(t_k)$. It follows that $d(1, n_k) \leq t_k$ for all $k$, so $n_k \to 1$ in $S$. Since $(S', d')$ is complete, it follows that the metric spaces $(S, d)$ and $(S', d')$ are isometric, except that the point labeled $1$ in $S$ is unlabeled in $S'$. Thus, we can also view $\gamma$ as a probability measure on $S'$.

Fix $t > 0$. Let $n_1, \ldots, n_{N(t)}$ be the smallest integers in the $N(t)$ blocks of $\Pi'(t)$. Note that if $x \in B(n_k, t)$, then $B(x, s) = B(n_k, s)$ for all $s \geq t$. Since $X$ has distribution $\gamma$, the probability that $X$ is in $B(n_k, s)$, conditional on $(\Pi'(t), t \geq 0)$, equals the asymptotic frequency of the block of $\Pi(t)$ containing $n_k$. Likewise, since $\Pi(t)$ is an exchangeable random partition, the probability that $1$ is in $B(n_k, s)$, conditional on $(\Pi'(t), t \geq 0)$, is the asymptotic frequency of the block of $\Pi(t)$ containing $n_k$. Since, whenever $x, y \in B(n_k, t)$, we have $\gamma(B(x,s)) = \gamma(B(y,s))$ for all $s \geq t$, it follows that the processes $(\gamma(B(1,s)), s \geq t)$ and $(\gamma(B(X,s)), s \geq t)$ have the same law. The lemma follows because $K(s) = \gamma(B(1,s))$ for all $s$. $\square$

The next proposition, combined with Proposition 6.1, proves Theorem 1.7.

**Proposition 6.4.** *We have $\dim_H(S) \geq 1/(\alpha - 1)$ a.s.*

**Proof.** Fix $\beta < 1/(\alpha - 1)$. We need to show that $\dim_H(S) \geq \beta$, and by Lemmas 6.2 and 6.3, it suffices to show that for some constant $C$, we have

$$\limsup_{t \downarrow 0} t^{-\beta} K(t) < C \quad \text{a.s.,} \tag{41}$$

where $K(t)$ is the asymptotic frequency of the block containing $1$ at time $t$.

Let $\Psi$ be a Poisson point process on $[0, \infty) \times \{0,1\}^\infty$ with intensity measure $dt \times L(d\xi)$, where $L(B) = \int_0^1 Q_x(B) x^{-2} \Lambda(dx)$ for all measurable $B$ and $Q_x$ is the distribution of an infinite sequence of $\{0,1\}$-valued random variables that are one with probability $x$ and zero with probability $1 - x$. We may assume that the coalescent process $(\Pi(t), t \geq 0)$ is constructed from $\Psi$ as described in Section 3. Choose a real number $b$ such that $\beta < b < 1/(\alpha - 1)$. Obtain a new Poisson point process $\Psi^*$ by removing from $\Psi$ all points $(t, \xi)$ such that $\xi_1 = 1$ and

$$\lim_{n \to \infty} \frac{1}{n} \sum_{i=1}^n \xi_i \geq t^b.$$

Then define a new coalescent process $(\Pi^*(t), t \geq 0)$ from $\Psi^*$, again using the procedure described in Section 3. Choose $\delta > 0$ such that $f(x) \leq 2Ax^{1-\alpha}$ for all $x \in (0, \delta]$. The number of points in $[0, t] \times \{0,1\}^\infty$ that are in $\Psi$ but not $\Psi^*$ has a Poisson distribution with mean

$$\int_0^t \int_{s^b}^1 x^{-1} \Lambda(dx) \, ds \leq \int_0^t \left( \delta^{-1} \Lambda([\delta, 1]) + \int_{s^b}^\delta 2Ax^{-\alpha} \, dx \right) ds$$

$$\leq \delta^{-1} \Lambda([\delta,1])t + \frac{2A}{(\alpha-1)(b(1-\alpha)+1)} t^{b(1-\alpha)+1}.$$

Since $b(1-\alpha) + 1 > 0$, this expression goes to zero as $t \downarrow 0$. It follows that there almost surely exists a $t^* > 0$ such that the processes $\Psi$ and $\Psi^*$ are the same on $[0, t^*] \times \{0, 1\}^\infty$. Therefore, $\Pi(t) = \Pi^*(t)$ for all $t \leq t^*$, so it suffices to show (41) with $K(t)$ replaced by $K^*(t)$, where $K^*(t)$ is the asymptotic frequency of the block containing 1 in $\Pi^*(t)$.

The partition $\Pi^*(t)$ is not exchangeable because only points $(t, \xi)$ with $\xi_1 = 1$ are removed from $\Psi$. However, the sequence whose $k$th term is the indicator of the event that 1 and $k+1$ are in the same block of $\Pi^*(t)$ is exchangeable. Therefore, the asymptotic frequency $K^*(t)$ of the block containing 1 exists almost surely, and its expected value is the probability that 1 and 2 are in the same block of $\Pi^*(t)$. This probability is bounded by the expected number of points $(s, \xi)$ of $\Psi^*$ such that $s \leq t$ and $\xi_1 = \xi_2 = 1$. Therefore, for $t$ small enough that $t^b \leq \delta$,

$$E[K^*(t)] \leq \int_0^t \int_0^{s^b} \Lambda(dx) \, ds \leq t\Lambda((0, t^b]) \leq 2At \int_0^{t^b} x^{1-\alpha} \, dx = \frac{2A}{2-\alpha} t^{1+b(2-\alpha)}.$$

Therefore, by Markov's Inequality, for all $\varepsilon > 0$, there exists a constant $C^*$ such that for all $t > 0$, we have $P(t^{-\beta} K^*(t) \geq \varepsilon) \leq C^* t^\eta$, where $\eta = 1 + b(2-\alpha) - \beta \geq 1 + b(2-\alpha) - b = 1 + b(1-\alpha) > 0$. It follows that

$$\sum_{k=1}^\infty P(2^{k\beta} K^*(2^{-k}) \geq \varepsilon) \leq C^* \sum_{k=1}^\infty 2^{-k\eta} < \infty.$$

By the Borel–Cantelli lemma, we have $2^{k\beta} K^*(2^{-k}) < \varepsilon$ for sufficiently large $k$ almost surely. Since the process $(K(t), t \geq 0)$ is nondecreasing, this implies (41) with $K^*$ in place of $K$. □

## 7. Dynamics of the number of blocks

In this section, we prove Theorems 1.8 and 1.9. Our first lemma gives a bound on the probabilities $\zeta_{n,k}$ which is uniform in $n$, under the additional assumption that the density of $\Lambda$ is bounded.

**Lemma 7.1.** *Assume that the assumptions of Theorem 1.8 hold, and that in addition there are constants $0 < C_1 \leq C_2 < \infty$ such that the function $f$ satisfies*

$$C_1 x^{1-\alpha} \leq f(x) \leq C_2 x^{1-\alpha} \tag{42}$$

*for all $x \in (0, 1]$. Then there exists a constant $C$ such that $\zeta_{n,k} \leq C k^{-1-\alpha}$ for all positive integers $n$ and $k$ such that $k \leq n - 1$.*

**Proof.** For $2 \leq k \leq n$, we have

$$\binom{n}{k} \lambda_{n,k} \leq C_2 \binom{n}{k} \int_0^1 x^{k-1-\alpha}(1-x)^{n-k} \, dx = \frac{C_2 n! \Gamma(k-\alpha)}{k! \Gamma(n-\alpha+1)} \leq C_3 n^\alpha k^{-1-\alpha}, \tag{43}$$

where $C_3$ is a constant that does not depend on $n$ or $k$. The same argument gives $\binom{n}{k} \lambda_{n,k} \geq C_4 n^\alpha k^{-1-\alpha}$, where $C_4$ is another constant, so

$$\lambda_n = \sum_{k=2}^n \binom{n}{k} \lambda_{n,k} \geq C_4 n^\alpha \sum_{k=2}^n k^{-1-\alpha} \geq C_5 n^\alpha \tag{44}$$

for some constant $C_5$. Eqs (43) and (44) give the result. □

**Lemma 7.2.** *Let $V_{n,k}$ be the event that $N(t) \in \{n, n+1, \ldots, n+k\}$ for some $t$. Under the assumptions of Theorem 1.8, for all $\varepsilon > 0$ there exists a positive integer $K$ such that for all positive integers $n$, we have $P(V_{n,n+k}^c) \leq \varepsilon$ whenever $k \geq K$.*

**Proof.** We first assume that the density $f$ satisfies (42) for all $x \in (0, 1]$. Recall that $V_n = V_{n,0}$ is the event that $N(t) = n$ for some $t$. For all positive integers $m$, let $B_m$ be a random variable such that $P(B_m = j) = \zeta_{m,j}$ for all $j \leq m - 1$. Note that

$$P(V_{n,n+k}^c \cap V_{n+k+1}) \leq P(V_{n,n+k}^c | V_{n+k+1}) = P(B_{n+k+1} > k+1) = \sum_{j=k+2}^{n+k} \zeta_{n+k+1,j} \leq \frac{C}{\alpha}(k+1)^{-\alpha}, \quad (45)$$

where $C$ is the constant from Lemma 7.1. Note that almost surely the coalescent must have more than $n$ blocks for sufficiently small $t$ because the coalescent process restricted to $\{1, \ldots, n+1\}$ almost surely holds in its initial state for a positive amount of time. Therefore, if $V_{n,n+k}^c$ occurs, then $V_{n,n+k+i}^c \cap V_{n+k+i+1}$ must occur for some nonnegative integer $i$. By (45),

$$P(V_{n,n+k}^c) \leq \sum_{i=0}^{\infty} P(V_{n,n+k+i}^c \cap V_{n+k+i+1}) \leq \frac{C}{\alpha} \sum_{i=0}^{\infty} (k+i+1)^{-\alpha} \leq \frac{C}{\alpha(\alpha-1)} k^{1-\alpha}.$$

Because this bound does not depend on $n$, the conclusion of the lemma holds when $f$ satisfies (7.1) for all $x \in (0, 1]$.

Next, consider the general case in which $f$ could be unbounded. The assumptions imply that there exists a $\delta > 0$ such that $f$ satisfies (42) for $x \in (0, \delta]$. Define a measure $\Lambda_1$ on $[0, 1]$ by $\Lambda_1(dx) = (f(x)\mathbf{1}_{\{x \leq \delta\}} + x^{1-\alpha}\mathbf{1}_{\{\delta < x < 1\}})dx$. Let $(\Pi_1(t), t \geq 0)$ be a $\Lambda_1$-coalescent. By Lemma 3.1, we may assume that the coalescent processes $(\Pi(t), t \geq 0)$ and $(\Pi_1(t), t \geq 0)$ are coupled so that there almost surely exists a random time $t > 0$ such that $\Pi(s) = \Pi_1(s)$ for all $s < t$. It follows that there exists a fixed time $u > 0$ such that $P(\Pi(s) = \Pi_1(s)$ for all $s \leq u) > 1 - \varepsilon/3$. Furthermore, there exists an integer $M$ such that $P(N(u) \leq M) \geq 1 - \varepsilon/3$. By our result for the case in which $f$ satisfies (42) for all $x \in (0, 1]$, there exists an integer $L$ such that whenever $k \geq L$, the probability that $\Pi_1(t) \in \{n, n+1, \ldots, n+k\}$ for some $t$ is at least $1 - \varepsilon/3$ for all $n$. However, if $\Pi(s) = \Pi_1(s)$ for all $s \leq u$ and $N(u) \leq k$, then $V_{n,n+k}$ occurs if and only if $\Pi_1(t) \in \{n, n+1, \ldots, n+k\}$ for some $t$. The result now follows by taking $K = \max\{L, M\}$. $\square$

**Proof of Theorem 1.8.** Let $(X_i)_{i=1}^{\infty}$ be a sequence of i.i.d. random variables with distribution $\zeta$. Let $S_0 = 0$, and for all positive integers $n$, let $S_n = \sum_{i=1}^{n} X_i$. For positive integers $k$, let $M_k = \max\{n \geq 0 : S_n \leq k\}$, and then define the age by $A_k = k - S_{M_k}$. The process $(A_k)_{k=0}^{\infty}$ is an irreducible Markov chain on the nonnegative integers, and the distribution of $A_k$ converges to a stationary distribution as $k \to \infty$. Since the distribution $\zeta$ has mean $1/(\alpha - 1)$, the expected time for the Markov chain to return to zero is $1/(\alpha - 1)$. It follows that

$$\lim_{k \to \infty} P(A_k = 0) = \alpha - 1.$$

Let $\varepsilon > 0$. Choose $K$ as in Lemma 7.2, so that $P(V_{n,n+k}) \geq 1 - \varepsilon$ whenever $k \geq K$. Choose an integer $m$ sufficiently large that whenever $k \geq (m-1)K$, we have $|P(A_k = 0) - (\alpha - 1)| \leq \varepsilon$.

Let $W_n = V_{n+(m-1)K, n+mK}$ be the event that $N(t) \in \{n+(m-1)K, \ldots, n+mK\}$ for some $t$. Let $H_n = \max\{j \leq n+mK : N(t) = j$ for some $t\}$. Note that $H_n \geq n + (m-1)K$ if and only if $W_n$ occurs. Let $(B_k)_{k=2}^{\infty}$ be a sequence of independent random variables such that $P(B_k = j) = \zeta_{k,j}$ for all positive integers $j \leq k - 1$. For nonnegative integers $i$ such that $i \leq K$ and positive integers $b_1, b_2, \ldots, b_{mK-i}$, define the function $F_i(b_1, b_2, \ldots, b_{mK-i})$ as follows. Construct a sequence $(x_k)_{k=0}^{mK-i}$ such that $x_0 = 0$ and, for $k \geq 0$, we have $x_{k+1} = x_k + b_{x_k}$ if $x_k \leq mK - i$ and $x_{k+1} = x_k$ otherwise. Then define $F_i(b_1, b_2, \ldots, b_{mK-i})$ to be 1 if $x_j = mK - i$ for some $j$, and 0 otherwise. Thinking of $B_k$ as the number of blocks lost in the next collision if the coalescent has $k$ blocks, we see that, conditional on the event that the coalescent has $n + mK - i$ blocks at some time, the probability that the coalescent eventually has exactly $n$ blocks is $E[F_i(B_{n+mK-i},\ldots,B_{n+1})]$. It follows that

$$P(V_n \cap W_n) = \sum_{i=0}^{K} P(H_n = n + mK - i) E[F_i(B_{n+mK-i},\ldots,B_{n+1})]. \tag{46}$$

Note also that

$$F_i(X_1, X_2, \ldots, X_{mK-i}) = \mathbf{1}_{\{A_{mK-i}=0\}}.$$

By (9), the distribution of $B_n$ converges to $\zeta$ as $n \to \infty$. Therefore, there is a constant $n_0$ such that for $n \geq n_0$, the sequences $(B_k)_{k=2}^{\infty}$ and $(X_i)_{i=1}^{\infty}$ can be coupled so that the total variation distance between $(B_{n+mK-i}, \ldots, B_{n+1})$ and $(X_1, \ldots, X_{mK-i})$ is at most $\varepsilon$ for all $i \leq K$. It follows that for $n \geq n_0$, we have

$$|E[F_i(B_{n+mK-i},\ldots,B_{n+1})] - P(A_{mK-i} = 0)| \leq \varepsilon.$$

Therefore by our choice of $m$, we have $|E[F_i(B_{n+mK-i},\ldots,B_{n+1})] - (\alpha - 1)| \leq 2\varepsilon$. Since $P(W_n) \geq 1 - \varepsilon$ by our choice of $K$, Eq. (46) now yields, for $n \geq n_0$,

$$P(V_n) \geq P(V_n \cap W_n) \geq P(W_n)(\alpha - 1 - 2\varepsilon) \geq (1 - \varepsilon)(\alpha - 1 - 2\varepsilon)$$

and

$$P(V_n) \leq P(W_n^c) + P(V_n \cap W_n) \leq \varepsilon + (\alpha - 1 + 2\varepsilon) = \alpha - 1 + 3\varepsilon.$$

The theorem follows by letting $\varepsilon \downarrow 0$. □

**Proof of Theorem 1.9.** Recall that $T_k = \inf\{t\colon N(t) \leq k\}$. Note that the distribution of $L_n$ is the same as the conditional distribution of

$$\int_{T_n}^{T_1} N(s)\,ds$$

given $N(T_n) = n$. Suppose $g\colon (0,\infty) \to (0,\infty)$ is a nonincreasing function such that $g(t) \sim Ct^{-\gamma}$ for some $\gamma > 1$, where $\sim$ means that the ratio of the two sides tends to zero as $t \downarrow 0$. Then it is straightforward to verify that for any $D > 0$, we have

$$\int_t^D g(s)\,ds \sim \frac{C}{\gamma - 1} t^{1-\gamma}. \tag{47}$$

By Theorem 1.1, we can apply (47) with $g(t) = N(t)$, $C = [\alpha/(A\Gamma(2-\alpha))]^{1/(\alpha-1)}$, and $\gamma = 1/(\alpha-1)$ to get

$$\int_t^{T_1} N(s)\,ds \sim \left(\frac{\alpha}{A\Gamma(2-\alpha)}\right)^{1/(\alpha-1)} \left(\frac{\alpha-1}{2-\alpha}\right) t^{(\alpha-2)/(\alpha-1)} \quad \text{a.s.} \tag{48}$$

Theorem 1.1 also implies that $T_n \sim [\alpha/(A\Gamma(2-\alpha))]n^{1-\alpha}$ a.s., where $\sim$ means that the ratio of the two sides tends to 1 as $n \to \infty$. Combining this observation with (48), we get

$$\lim_{n\to\infty} \frac{1}{n^{2-\alpha}} \int_{T_n}^{T_1} N(s)\,ds = \frac{\alpha(\alpha-1)}{A\Gamma(2-\alpha)(2-\alpha)} \quad \text{a.s.} \tag{49}$$

Let $\varepsilon > 0$. By (49) and Theorem 1.8, there exists an $M$ such that if $n \geq M$, then $P(V_n) = P(N(T_n) = n) \geq (\alpha - 1)/2$ and

$$P\left(\left|\frac{1}{n^{2-\alpha}} \int_{T_n}^{T_1} N(s)\,ds - \frac{\alpha(\alpha-1)}{A\Gamma(2-\alpha)(2-\alpha)}\right| > \varepsilon\right) \leq \frac{\varepsilon(\alpha-1)}{2}.$$

Since $L_n$ is the conditional distribution of $\int_{T_n}^{T_1} N(s)\,\mathrm{d}s$ given $V_n$, for $n \geq M$ we have

$$P\left(\left|\frac{L_n}{n^{2-\alpha}} - \frac{\alpha(\alpha-1)}{A\Gamma(2-\alpha)(2-\alpha)}\right| > \varepsilon\right) \leq \frac{\varepsilon(\alpha-1)}{2P(V_n)} \leq \varepsilon.$$

The result follows. □

## Acknowledgments

The authors thank Jean Bertoin and Anja Sturm for helpful discussions and a referee for helpful comments. N.B. thanks the L.A.T.P. in Marseille and U.C.S.D. for their invitations while this work was being done.

## References

[1] D. J. Aldous. Deterministic and stochastic models for coalescence (aggregation and coagulation): a review of the mean-field theory for probabilists. *Bernoulli* **5** (1999) 3–48. MR1673235

[2] A.-L. Basdevant. Ruelle's probability cascades seen as a fragmentation process. *Markov Process. Related Fields* **12** (2006) 447–474. MR2246260

[3] J. Berestycki, N. Berestycki and J. Schweinsberg. Beta-coalescents and continuous stable random trees. *Ann. Probab.* **35** (2007) 1835–1887.

[4] J. Bertoin. *Lévy Processes*. Cambridge University Press, Cambridge, 1996. MR1406564

[5] J. Bertoin. *Random Coagulation and Fragmentation Processes*. Cambridge University Press, Cambridge, 2006. MR2253162

[6] J. Bertoin and J.-F. Le Gall. The Bolthausen–Sznitman coalescent and the genealogy of continuous-state branching processes. *Probab. Theory Related Fields* **117** (2000) 249–266. MR1771663

[7] J. Bertoin and J.-F. Le Gall. Stochastic flows associated to coalescent processes. *Probab. Theory Related Fields* **126** (2003) 261–288. MR1990057

[8] J. Bertoin and J.-F. Le Gall. Stochastic flows associated to coalescent processes III: limit theorems. *Illinois J. Math.* **50** (2006) 147–181. MR2247827

[9] J. Bertoin and J. Pitman. Two coalescents derived from the ranges of stable subordinators. *Electron. J. Probab.* **5** (2000) 1–17. MR1768841

[10] N. H. Bingham, C. M. Goldie and J. L. Teugels. *Regular Variation*. Cambridge University Press, Cambridge, 1987. MR0898871

[11] M. Birkner, J. Blath, M. Capaldo, A. Etheridge, M. Möhle, J. Schweinsberg and A. Wakolbinger. Alpha-stable branching and beta-coalescents. *Electron. J. Probab.* **10** (2005) 303–325. MR2120246

[12] E. Bolthausen and A.-S. Sznitman. On Ruelle's probability cascades and an abstract cavity method. *Comm. Math. Phys.* **197** (1998) 247–276. MR1652734

[13] P. Donnelly, S. N. Evans, K. Fleischmann, T. G. Kurtz and X. Zhou. Continuum-sites stepping-stone models, coalescing exchangeable partitions, and random trees. *Ann. Probab.* **28** (2000) 1063–1110. MR1797304

[14] P. Donnelly and T. G. Kurtz. Particle representations for measure-valued population models. *Ann. Probab.* **27** (1999) 166–205. MR1681126

[15] R. M. Dudley. *Real Analysis and Probability*. Wadsworth and Brooks/Cole, Pacific Grove, CA, 1989. MR0982264

[16] T. Duquesne. A limit theorem for the contour process of conditioned Galton–Watson trees. *Ann. Probab.* **31** (2003) 996–1027. MR1964956

[17] T. Duquesne and J.-F. Le Gall. Probabilistic and fractal aspects of Lévy trees. *Probab. Theory Related Fields* **131** (2005) 553–603. MR2147221

[18] R. Durrett and J. Schweinsberg. A coalescent model for the effect of advantageous mutations on the genealogy of a population. *Stochastic Process. Appl.* **115** (2005) 1628–1657. MR2165337

[19] N. El Karoui and S. Roelly. Propriétés de martingales, explosion et représentation de Lévy–Khintchine d'une classe de processus de branchement à valeurs mesures. *Stochastic Process. Appl.* **38** (1991) 239–266. MR1119983

[20] S. N. Evans. Kingman's coalescent as a random metric space. In *Stochastic Models: A Conference in Honour of Professor Donald A. Dawson* (L. G. Gorostiza and B. G. Ivanoff, Eds). Canadian Mathematical Society/American Mathematical Society, 2000. MR1775475

[21] K. Falconer. *Fractal Geometry: Mathematical Foundations and Applications*, 2nd edition. Wiley, Hoboken, NJ, 2003. MR2118797

[22] C. Goldschmidt and J. Martin. Random recursive trees and the Bolthausen–Sznitman coalescent. *Electron. J. Probab.* **10** (2005) 718–745. MR2164028

[23] J. Kesten, P. Ney and F. Spitzer. The Galton–Watson process with mean one and finite variance. *Theory Probab. Appl.* **11** (1966) 513–540. MR0207052

[24] J. F. C. Kingman. The representation of partition structures. *J. London Math. Soc.* **18** (1978) 374–380. MR0509954

[25] J. F. C. Kingman. The coalescent. *Stochastic Process. Appl.* **13** (1982) 235–248. MR0671034
[26] J. Lamperti. The limit of a sequence of branching processes. *Z. Wahrsch. Verw. Gebiete* **7** (1967) 271–288. MR0217893
[27] J. Lamperti. Continuous state branching processes. *Bull. Amer. Math. Soc.* **73** (1967) 382–386. MR0208685
[28] M. Möhle. On the number of segregating sites for populations with large family sizes. *Adv. in Appl. Probab.* **38** (2006) 750–767.
[29] M. Möhle and S. Sagitov. A classification of coalescent processes for haploid exchangeable population models. *Ann. Probab.* **29** (2001) 1547–1562. MR1880231
[30] J. Pitman. Coalescents with multiple collisions. *Ann. Probab.* **27** (1999) 1870–1902. MR1742892
[31] S. Sagitov. The general coalescent with asynchronous mergers of ancestral lines. *J. Appl. Probab.* **36** (1999) 1116–1125. MR1742154
[32] J. Schweinsberg. A necessary and sufficient condition for the $\Lambda$-coalescent to come down from infinity. *Electron. Comm. Probab.* **5** (2000) 1–11. MR1736720
[33] J. Schweinsberg. Coalescent processes obtained from supercritical Galton–Watson processes. *Stochastic Process. Appl.* **106** (2003) 107–139. MR1983046
[34] M. L. Silverstein. A new approach to local times. *J. Math. Mech.* **17** (1968) 1023–1054. MR0226734
[35] R. Slack. A branching process with mean one and possibly infinite variance. *Z. Wahrsch. Verw. Gebiete* **9** (1968) 139–145. MR0228077
[36] A. M. Yaglom. Certain limit theorems of the theory of branching processes. *Dokl. Acad. Nauk SSSR* **56** (1947) 795–798. MR0022045